\documentclass[11pt]{article}
\usepackage{amsfonts}
\usepackage{amsfonts}
\usepackage{amsfonts}
\usepackage{mathtools}
\usepackage{mathrsfs}
\usepackage{bm}
\usepackage{amssymb,amsmath} 
 \allowdisplaybreaks

\catcode`!=11
\let\!int\int \def\int{\displaystyle\!int}
\let\!lim\lim \def\lim{\displaystyle\!lim}
\let\!sum\sum \def\sum{\displaystyle\!sum}
\let\!sup\sup \def\sup{\displaystyle\!sup}
\let\!inf\inf \def\inf{\displaystyle\!inf}
\let\!cap\cap \def\cap{\displaystyle\!cap}
\let\!max\max \def\max{\displaystyle\!max}
\let\!min\min \def\min{\displaystyle\!min}
\let\!frac\frac \def\frac{\displaystyle\!frac}
\catcode`!=12

\let\oldsection\section
\renewcommand\section{\setcounter{equation}{0}\oldsection}

\allowdisplaybreaks
\def\pf{\it{Proof.}\rm\quad}
\newcommand\kk{\mbox{\bfseries{k}}}
\newcommand\BB{\mbox{\bfseries{B}}}
\newcommand\vv{\mbox{\bfseries{v}}}\newcommand\ff{\mbox{\bfseries{f}}}
\newcommand\uu{\mbox{\bfseries{u}}}

\newcommand\divg{{\text{div}}}

\newtheorem{thm}{Theorem}[section]
\newtheorem{pro}{Proposition}[section]

\begin{document}
\title {\bf On the infinite Prandtl number limit in two-dimensional magneto-convection }
\author{
{Jianwen Zhang\thanks{Email: jwzhang@xmu.edu.cn}}\quad {and}\quad {Mingyu Zhang\thanks{Corresponding author. Email: myuzhang@126.com}}\\[1mm]
\small School of Mathematical Sciences, Xiamen University,  Xiamen
361005, P.R. China}

\date{}
\maketitle \noindent{\bf Abstract.} In this paper, the infinite  limit of the Prandtl number is justified for the two-dimensional incompressible magneto-convection, which describes the nonlinear interaction between the
Rayleigh-B$\rm\acute{e}$nard convection and an externally magnetic field. Both the convergence rates and the thickness of initial layer are obtained. Moreover, based on the method of formal asymptotic expansions, an effective dynamics is constructed to simulate the motion within the initial layer.

\section{Introduction}
In this paper, we are interested in a two-dimensional Boussinesq fluid with nonlinear interaction between
Rayleigh-B$\rm\acute{e}$nard convection and an externally magnetic field. To begin, let us consider a
horizontally stratified fluid layer of {\it characteristic height}
$h$, referred to as a Cartesian coordinate system with $x$-axis
in the horizontal direction and $y$-axis pointing vertically
upward. Assume that a fixed temperature difference, say
$\theta_2-\theta_1$, is maintained across the layer of the fluid heated from
below in an externally imposed magnetic field $\BB_0 =\bar B\kk$. For
simplicity, we also assume the periodic boundary conditions in the
horizontal direction. Then, the  MHD-Boussinesq
approximation for incompressible viscous and resistive flows
 reads as follows (cf. \cite{Da2001}):
\begin{equation}\label{1.1}
\begin{dcases}
\frac{\partial \uu}{\partial t}+\uu\cdot \nabla\uu+\nabla p=\mu\Delta\uu+g\alpha\kk \theta+\BB\cdot\nabla\BB,\\
\frac{\partial \BB}{\partial t}+\uu\cdot \nabla\BB-\BB\cdot \nabla \uu=\nu\Delta \BB,\\
 \frac{\partial \theta}{\partial t}+\uu\cdot \nabla \theta=\kappa\Delta \theta,\\
\nabla\cdot\uu=0,\quad\nabla\cdot\BB=0
\end{dcases}
\end{equation}
with the following initial and boundary conditions:
\begin{equation}\label{1.2}
\begin{dcases}
\uu|_{t=0}=\uu_0,\quad \BB|_{t=0}=\BB_0,\quad \theta|_{t=0}=\theta_0,\\
\uu|_{y=0,h}=0,\quad \BB|_{y=0,h}=\bar B \kk,\quad \theta|_{y=0}=\theta_2,\quad \theta|_{y=h}=\theta_1,
\end{dcases}
\end{equation}
and the periodic boundary conditions in the horizontal direction.
Here, $\uu=(u_1,u_2)$ is the velocity, $\BB$ is the magnetic field, $\theta$ is the temperature and $p$ is the total pressure (incorporating the magnetic pressure); $\mu$ is the  kinematic viscosity, $\nu$ is the magnetic diffusivity, $\kappa$ is the thermal diffusivity, $g$ is the gravitational acceleration, $\alpha$ is the thermal expansion coefficient, and $\kk$ is the vertical unit vector.

This set of equations describes the nonlinear  interaction between Rayleigh-B${\rm\acute{e}}$nard convection and an externally imposed magnetic fields, which is called magneto-convection and may explain certain prominent features on the solar surface. It was shown in \cite{Bu1989,G1998} that the finite amplitude onset of steady convection became possible when the Rayleigh number is considerably below the values predicted by linear theory. Magnetic fields with sunspots are sufficiently strong to suppress convection on granular and supergranular scales (see \cite{C1952,C1954,C1961,N1955,N1957,Pa1979,T1951}). However, we are far from a real understanding of the dynamic coupling between convection and magnetic fields in stars and magnetically confined high-temperature plasmas. So, it is of great importance to understand how energy  transport and convection are affected by an imposed magnetic field, that is, how the Lorentz force affects the convection patterns in sunspots and magnetically confined high-temperature plasmas.

Since we aim to consider the problem of magneto-convection, it is more convenient to consider the standard and natural non-dimensional equations of (\ref{1.1}). So, if using the units of the layer height $h$ as the characteristic length scale, the thermal diffusion time $h^2/\kappa$ as the characteristic time scale, the ratio of typical length over typical time $\kappa/h$ as the typical velocity, the imposed field strength $\bar B$ as the typical magnetic field, and the temperature on a scale where the upside is kept at $0$ and the downside is kept at $1$, then we obtain the  following non-dimensional equations of (\ref{1.1}):
\begin{equation}\label{1.3}
\begin{dcases}
\frac{1}{\rm Pr}\left(\frac{\partial \uu}{\partial t}+\uu\cdot \nabla\uu\right)+\nabla p=\Delta\uu+{\rm Ra}\kk \theta+{\rm Pm\cdot Q}\BB\cdot\nabla\BB+{\rm Pm \cdot Q}\frac{\partial\BB}{\partial y},\\
\frac{1}{\rm Pm}\left(\frac{\partial \BB}{\partial t}+\uu\cdot \nabla\BB-\BB\cdot \nabla \uu-\frac{\partial \uu}{\partial y}\right)=\Delta \BB,\\
 \frac{\partial \theta}{\partial t}+\uu\cdot \nabla \theta=\Delta \theta,\\
\nabla\cdot\uu=0,\quad\nabla\cdot\BB=0
\end{dcases}
\end{equation}
with periodic boundary conditions in the horizontal direction and 
\begin{equation}\label{1.4}
\begin{dcases}
\uu|_{t=0}=\uu_0,\quad \BB|_{t=0}=\BB_0,\quad \theta|_{t=0}=\theta_0,\\
\uu|_{y=0,1}=0,\quad \BB|_{y=0,1}=0,\quad \theta|_{y=0}=1,\quad \theta|_{y=1}=0.
\end{dcases}
\end{equation}
Here, for simplicity but without any confusion, we still use the notations $(\uu,\BB,\theta)$ and $(\uu_0,\BB_0,\theta_0)$ to denote the non-dimensional quantities and the initial data, respectively.

There are  four important dimensionless parameters in (\ref{1.3}): the {\it Rayleigh number}
$$
{\rm Ra}=\frac{g\alpha (\theta_2-\theta_1)h^3}{\mu\kappa},
$$
measuring the ratio of overall buoyancy force to the damping coefficients; the {\it Chandrasekhar number}
$$
{\rm Q}=\frac{\bar B^2h^2}{\mu\nu},
$$
measuring the ratio of Lorentz force to viscosity; the {\it Prandtl number}
$$
{\rm Pr}=\frac{\mu}{\kappa},
$$
measuring the relative ratio of momentum diffusivity to thermal diffusivity; and the {\it magnetic Prandtl number}
$$
{\rm Pm}=\frac{\nu}{\kappa},
$$
measuring the ratio of magnetic diffusivity to thermal diffusivity.

The problem of large Prandtl number (i.e., ${\rm Pr}\gg 1$) finds its important many applications for the fluids such as silicone oil, the earth's mentls, as well as the gases under high pressure (see, for example, \cite{Bu1989,C1961,GL2000,S1994}). Since we have normalized the time to the thermal diffusive time scale, the large Prandtl number means that the viscous time scale of the fluid (i.e., $h^2/\mu$) is much shorter than the thermal diffusive time scale (i.e., $h^2/\kappa$). Thus, the velocity field slaved by the temperature field will settle into some ``equilibrium" state due to the long-time viscosity effect (cf. \cite{Wang2004,Wang2007}). Formally, if the Prandtl number equal to infinity in (\ref{1.3}), then the convection term can be negligible and the the so-called {\it infinite Prandtl number system} reads 
\begin{equation}\label{1.5}
\begin{dcases}
\nabla p^0-\Delta\uu^0={\rm Ra}\kk \theta^0+{\rm Pm\cdot Q}\BB^0\cdot\nabla\BB^0+{\rm Pm\cdot Q}\frac{\partial\BB^0}{\partial y},\\
\frac{1}{\rm Pm}\left(\frac{\partial \BB^0}{\partial t}+\uu^0\cdot \nabla\BB^0-\BB^0\cdot \nabla \uu^0-\frac{\partial \uu^0}{\partial y}\right)=\Delta \BB^0,\\
 \frac{\partial \theta^0}{\partial t}+\uu^0\cdot \nabla \theta^0=\Delta \theta^0,\\
\nabla\cdot\uu^0=0,\quad\nabla\cdot\BB^0=0,
\end{dcases}
\end{equation}
with periodic boundary conditions in the $x$-direction and
\begin{equation}\label{1.6}
\begin{dcases}
\BB^0|_{t=0}=\BB_0,\quad \theta^0|_{t=0}=\theta_0,\\
\uu^0|_{y=0,1}=0,\quad \BB^0|_{y=0,1}=0,\quad \theta^0|_{y=0}=1,\quad \theta^0|_{y=1}=0.
\end{dcases}
\end{equation}

The large Prandtl number for the incompressible fluids without magnetic effects has been studied by many authors (see, for example, \cite{CD1999,CHP2001,DC2001,Wang2004,Wang2007}). Indeed, if we ignore the magnetic effects in (\ref{1.5}), then  it turns into
 \begin{equation}\label{1.7}
\begin{dcases}
\nabla p^0-\Delta\uu^0={\rm Ra}\kk \theta^0,\\
 \frac{\partial \theta^0}{\partial t}+\uu^0\cdot \nabla \theta^0=\Delta \theta^0,\\
\nabla\cdot\uu^0=0.
\end{dcases}
\end{equation}
It is easy to utilize the standard results of the Stokes equations (cf. \cite{L1963,T2001}) to show that there exists a global strong solution to the initial and boundary value problem of (\ref{1.7}) even in the three-dimensional setting, since it readily follows from the maximum principle that the temperature $\theta^0$ is globally bounded. However, there is very few  result about the large Prandtl number for magneto-convection. In fact, it is easy to see that system (\ref{1.5})  is a coupled parabolic-elliptic system, which retains the essentially nonlinear Lorentz effect on the fluid. So, the rigorous mathematical theory (e.g., well-posedness, asymptotic behavior, {\it etc}) of (\ref{1.5}) in the three-dimensional framework is full of challenge,  though the stabilizing effect of magnetic field has been exploited in many works both from the physical and from the numerical point of view (see, for example, \cite{Ga1985,MR1994,MR2003,R1997}).

The global well-posedness theory of the equations for incompressible viscous fluids is classical and well-known, see, for example, \cite{L1963,ST1983,T2001} and the references cited therein. Moreover, a similar system as that in (\ref{1.5}) was also studied in \cite{Br2014,MRR2014} and the global existence of weak solutions was proved. We also mention the interesting works \cite{CW2011,CWY2014}, where the two-dimensional incompressible MHD equations with partial viscosities were considered. The main purpose of this paper is to justify the global-in-time asymptotic limit from the two-dimensional system (\ref{1.3}) to the one (\ref{1.5}) rigorously, as the Prandtl number ${\rm Pr}$ tends to infinity (i.e., $\Pr\to\infty$). As a by-product, the global well-posedness of strong solutions to the problem (\ref{1.5})--(\ref{1.6}) with large data was also proved. It is clear that the infinite limit of the Prandtl number is a singular one involving an initial layer.

Our first result is concerned with the convergence from $(\uu,\BB,\theta)$ to $(\uu^0,\BB^0,\theta^0)$ strictly away from the initial layer. Note that, since we are only interested in the infinite limit of Prandtl number,  for simplicity we assume   throughout the remainder of this paper that ${\rm Ra}, {\rm Q}, {\rm Pm}\equiv 1$.
\begin{thm}\label{thm1.1} Let $\Omega\triangleq [0,L]\times[0,1]$ and  ${\rm Ra, Q, Pm}\equiv1$.
Assume that
\begin{equation}\label{1.8}
\begin{dcases}
(\mathbf{u}_0, \mathbf{B}_0,\theta_0)\in H^2,\quad \nabla\cdot\mathbf{u}_0=0,\quad \nabla\cdot\mathbf{B}_0=0,\\
 (\mathbf{u}_0,\mathbf{B}_{0})|_{y=0,1}=0,\quad \theta_0|_{y=0}=1,\quad  \theta_0|_{y=1}=0.
 \end{dcases}
\end{equation}
Then for any  $0<T<\infty$, there exists a global unique solution $({\mathbf u},{\mathbf B},\theta)$ (resp. $({\mathbf u}^0,{\mathbf B^0},\theta^0)$) to the problem \eqref{1.3}--\eqref{1.4} (resp. the problem \eqref{1.5}--\eqref{1.6}) on $\Omega\times[0,T]$, such that
\begin{equation}\label{1.9}
 \sup_{0\leq t\leq T}\left(\|({\mathbf B}-{\mathbf B^0}, \theta-\theta^0)(t)\|^2_{H^1}+\|(\nabla({\mathbf B}-{\mathbf B^0}),\nabla(\theta-\theta^0))(t)\|^4_{L^4}\right)\leq C\varepsilon
\end{equation}
and
\begin{equation}\label{1.10}
\|(\mathbf{u}-\mathbf{u}^0)(t)\|_{H^2}^2\leq \frac{C\varepsilon}{t}\quad {for\ \ any}\quad \varepsilon^{1-\alpha}\leq t\leq T\ \ {with}\ \ \alpha\in(0,1),
\end{equation}
where $\varepsilon\triangleq {\rm Pr}^{-1}\in (0,1)$ and $C$ is a positive constant independent of $\varepsilon$.
\end{thm}

\noindent{\bf Remark 1.1} It seems unsatisfactory that the quantity $\|\BB\|_{H^2}$ cannot be uniformly bounded, although $\BB_0\in H^2$, and consequently, the convergence of $\BB-\BB^0$ in $H^2$ cannot be obtained. This is mainly due to the effects caused by the initial layer and the imposed strength  $\bar B \kk$ of
magnetic field, the latter of which induces an additional term $\BB_y$ in (\ref{1.3})$_2$.

\vskip 2mm

It is easily seen from (\ref{1.10})  that there is an initial layer between $\uu$ and $\uu^0$, whose thickness is almost of the value $\varepsilon$. Motivated by this fact, to capture the effective dynamics of the initial layer, we adopt the so-called two-time scale approach (cf. \cite{GL2000,Ho1995,Ma2003}) by introducing the following fast time scale
\begin{equation*}
\tau= {\rm Pr} \cdot t=\frac{t}{\varepsilon}\quad {\rm with} \quad \varepsilon=\frac{1}{\rm Pr}
\end{equation*}
and the formal asymptotic expansions
\begin{equation}\label{1.11}
\begin{dcases}
\uu=\uu^{(0)}(t,\tau)+\varepsilon \uu^{(1)}(t,\tau)+{\rm h.o.t.},\\
\BB(t,x)=\BB^{(0)}(t,\tau)+\varepsilon\BB^{(1)}(t,\tau)+{\rm h.o.t.},\\
\theta(t,x)=\theta^{(0)}(t,\tau)+\varepsilon \theta^{(1)}(t,\tau)+{\rm h.o.t.},
\end{dcases}
\end{equation}
where ``${\rm h.o.t.}$" represents the higher-order terms in $\varepsilon$.  Moreover, to ensure the validity of the formal asymptotic expansion for large values of the fast variable $\tau$, we also impose the customary sublinear growth condition
\begin{equation}\label{1.12}
\lim_{\tau\to \infty}\frac{(\uu^{(1)},\BB^{(1)},\theta^{(1)})(t,\tau)}{\tau}=0.
\end{equation}

Inserting the formal asymptotic expansion (\ref{1.11}) into  (\ref{1.3}) and noting that
$$
D_t=\frac{\partial}{\partial t}+\frac{1}{\varepsilon}\frac{\partial}{\partial\tau},
$$
we obtain after collecting the leading-order terms that (${\rm Ra, Q, Pm}\equiv1$)
\begin{equation}\label{1.13}
\begin{dcases}
\frac{\partial \uu^{(0)}}{\partial \tau}+\nabla p^{(0)}=\Delta \uu^{(0)}+ \kk \theta^{(0)}+\BB^{(0)}\cdot \nabla \BB^{(0)}+\frac{\partial{\BB^{(0)}}}{\partial y},\\
\frac{\partial\BB^{(0)}}{\partial \tau}=0,\quad
\frac{\partial \theta^{(0)}}{\partial \tau}=0,\quad\nabla\cdot\uu^{(0)}=0,\quad\nabla\cdot\BB^{(0)}=0.
\end{dcases}
\end{equation}

Let $A$ be the Stokes operator defined as
\begin{equation}\label{1.14}
A\uu={\mathbf f},
\end{equation}
if and only if $\uu$ satisfies
\begin{equation*}
\nabla p-\Delta \uu={\mathbf f},\quad\nabla\cdot\uu=0,\quad \uu|_{y=0,1}=0,
\end{equation*}
and the periodic conditions in the $x$-direction.

Then it follows from (\ref{1.13}) that
\begin{equation}\label{1.15}
\BB^{(0)}(t,\tau)=\BB^{(0)}(t),\quad \theta^{(0)}(t,\tau)=\theta^{(0)}(t),
\end{equation}
and
\begin{equation}\label{1.16}
\begin{split}
\uu^{(0)}(t,\tau)=&e^{-\tau A}\uu^{(0)}(t,0)+A^{-1}\mathbb{P}\left(\kk \theta^{(0)}+\BB^{(0)}\cdot\nabla \BB^{(0)}+\frac{\partial\BB^{(0)}}{\partial y}\right)(t)\\
\quad &-e^{-\tau A}A^{-1}\mathbb{P}\left(\kk \theta^{(0)}+\BB^{(0)}\cdot \nabla \BB^{(0)}+\frac{\partial\BB^{(0)}}{\partial y}\right)(t),
\end{split}
\end{equation}
where $\mathbb{P}$ is the Leray-Hopf projector and $A=-\mathbb{P}\Delta$ (see, for example, \cite{T2001}).

The next-order dynamics is governed by
\begin{equation}\label{1.17}
\begin{dcases}
\frac{\partial \uu^{(1)}}{\partial \tau}+A \uu^{(1)}=-\frac{\partial \mathbf{u}^{(0)}}{\partial t}-\mathbb{P}\left(\uu^{(0)}\cdot \nabla\uu^{(0)}\right)\\
\qquad\qquad\qquad\qquad+\mathbb{P}\left(\kk\theta^{(1)}+\BB^{(0)}\cdot
\nabla\BB^{(1)}+\BB^{(1)}\cdot \nabla\BB^{(0)}+\frac{\partial \BB^{(1)}}{\partial y}\right),\\
\frac{\partial\BB^{(1)}}{\partial \tau}=\Delta \BB^{(0)}-\uu^{(0)}\cdot \nabla \BB^{(0)}+\BB^{(0)}\cdot \nabla\uu^{(0)}+\frac{\partial \uu^{(0)}}{\partial y}-\frac{\partial \BB^{(0)}}{\partial t},\\
\frac{\partial \theta^{(1)}}{\partial \tau}=\Delta \theta^{(0)}-\uu^{(0)}\cdot \nabla \theta^{(0)}-\frac{\partial \theta^{(0)}}{\partial t}.
\end{dcases}
\end{equation}

In view of the sublinear growth condition (\ref{1.12}), we have
\begin{equation}\label{1.18}
\begin{dcases}
0=\Delta \BB^{(0)}-\uu^{(0)}\cdot \nabla \BB^{(0)}+\BB^{(0)}\cdot \nabla \uu^{(0)}+\frac{\partial \uu^{(0)}}{\partial y}-\frac{\partial \BB^{(0)}}{\partial t},\\
0=\Delta \theta^{(0)}-\uu^{(0)}\cdot \nabla \theta^{(0)}-\frac{\partial\theta^{(0)}}{\partial t},
\end{dcases}
\end{equation}
which is the limit model of infinite Prandtl number .

The equation (\ref{1.17})$_1$ for $\uu^{(1)}$ is dissipative. However, similarly to that in \cite{Wang2004}, there are three terms of $\uu^{(0)}$ in (\ref{1.16}): one term slaved by  the leading-order terms of temperature and Lorentz force, and another two terms exponentially decaying in time (initial layer type). This means that no more dynamics on $\uu^{(0)}$ is necessary except the ones in (\ref{1.16}). Moreover, by modifying the initial layer terms in such a way so that the initial data are fixed, we can propose the following {\it effective dynamics} within the initial layer ($\tau=t/\varepsilon$):
\begin{equation}\label{1.19}
\begin{dcases}
\uu^{(0)}=e^{-\tau A}\uu_0-e^{-\tau A}A^{-1}\mathbb{P}\left(\kk \theta_{0}+\BB_0\cdot \nabla \BB_0+\frac{\partial\BB_0}{\partial y}\right)\\
\qquad\quad+A^{-1}\mathbb{P}\left(\kk \theta^{(0)}+\BB^{(0)}\cdot\nabla \BB^{(0)}+\frac{\partial\BB^{(0)}}{\partial y}\right),\\
\frac{\partial \BB^{(0)}}{\partial t}+\uu^{(0)}\cdot \nabla \BB^{(0)}-\BB^{(0)}\cdot \nabla \uu^{(0)}-\frac{\partial \uu^{(0)}}{\partial y}=\Delta \BB^{(0)},\quad\nabla\cdot\BB^{(0)}=0,\\
\frac{\partial\theta^{(0)}}{\partial t}+\uu^{(0)}\cdot \nabla \theta^{(0)}=\Delta \theta^{(0)},
\end{dcases}
\end{equation}
which is completed with periodic boundary conditions in the $x$-direction and
\begin{equation}\label{1.20}
\begin{dcases}
\BB^{(0)}|_{t=0}=\BB_0,\quad \theta^{(0)}|_{t=0}=\theta_0,\\
\BB^{(0)}|_{y=0,1}=0,\quad\theta^{(0)}|_{y=0}=1,\quad\theta^{(0)}|_{y=1}=0.
\end{dcases}
\end{equation}
The solutions of (\ref{1.19})--(\ref{1.20}) will be compared with the ones of the problems (\ref{1.3})--(\ref{1.4}) and (\ref{1.5})--(\ref{1.6}) with ${\rm Ra, Q, Pm}\equiv1$, respectively.

With the help of the {\it effective dynamics} (\ref{1.19})--(\ref{1.20}), we can prove the following main theorem which is concerned with the behavior of the initial layer.
\begin{thm}\label{thm1.2} Let the conditions of Theorem \ref{thm1.1} be in force. For any fixed $0<T<\infty$, assume that
 $({\mathbf u},{\mathbf B},\theta)$,  $({\mathbf u}^0,{\mathbf B^0},\theta^0)$ and $({\mathbf u}^{(0)},{\mathbf B^{(0)}},\theta^{(0)})$ are the solutions of the problems \eqref{1.3}--\eqref{1.4},  \eqref{1.5}--\eqref{1.6} and \eqref{1.19}--\eqref{1.20} on $\Omega\times[0,T]$, respectively. Then, in addition to \eqref{1.9} and \eqref{1.10}, there exists a positive constant $C$, independent of $\varepsilon$, such that for any $t\in[0,T]$,
\begin{equation}\label{1.21}
\left\|\mathbf{B}-\mathbf{B}^0\right\|_{L^2}+\left\|\theta-\theta^0\right\|_{L^2}\leq C\varepsilon,
\end{equation}
\begin{equation}\label{1.22}
\left\|\mathbf{u}-\mathbf{u}^0-e^{-\tau A}\mathbf{u}_0+e^{-\tau A} A^{-1} \mathbb{P}\left(\mathbf{k} \theta_0+\mathbf{B}_0\cdot \nabla \mathbf{B}_0+\partial y \mathbf{B}_0\right)\right\|_{L^2}\leq  C\varepsilon
\end{equation}
and
\begin{equation}\label{1.23}
\left\|\nabla\left(\mathbf{u}-\mathbf{u}^0-e^{-\tau A}\mathbf{u}_0+e^{-\tau A} A^{-1} \mathbb{P}(\mathbf{k} \theta_0+\mathbf{B}_0\cdot \nabla \mathbf{B}_0+\partial y \mathbf{B}_0)\right)\right\|_{H^1}\leq  C\varepsilon^{1/2}.
\end{equation}
\end{thm}

\noindent{\bf Remark 1.2} It is easily seen from (\ref{1.22}) and (\ref{1.23}) that the motion within the initial layer can be modelled by the initial-layer correction function $e^{-\tau A}\mathbf{u}_0-e^{-\tau A} A^{-1} \mathbb{P}(\mathbf{k} \theta_0+\mathbf{B}_0\cdot \nabla \mathbf{B}_0+\partial y \mathbf{B}_0)$ in the sense of uniform convergence. Indeed, the first term $\mathbf{u}_0$ is the initial data of $\uu$ and the second one $A^{-1} \mathbb{P}(\mathbf{k} \theta_0+\mathbf{B}_0\cdot \nabla \mathbf{B}_0+\partial y \mathbf{B}_0)$ is the initial data of $\uu^0$.

\vskip 2mm

\noindent{\bf Remark 1.3} The $L^2$-convergence rates of order $\varepsilon$ in (\ref{1.21}) and (\ref{1.22}) are optimal, which can be justfied via the systematic asymptotic expansion with the small parameter $\varepsilon=1/{\rm Pr}$.

\vskip 2mm

The proofs of Theorems \ref{thm1.1} and \ref{thm1.2} are based on the global (uniform) estimates of  $({\mathbf u},{\mathbf B},\theta)$,  $({\mathbf u}^0,{\mathbf B^0},\theta^0)$ and $({\mathbf u}^{(0)},{\mathbf B^{(0)}},\theta^{(0)})$, which will be achieved by making a full use of the estimates of the Stokes equations in a manner similar to that used for the standard incompressible Navier-Stokes/MHD equations (cf. \cite{L1963,ST1983,T2001}). It is worth pointing out that due to the presence of initial layer, the uniform $L^p$-estimate of $\uu_t$ with any $p\geq1$ cannot be expected. In other words, the $L^p$-estimate of $\uu_t$ is actually $\varepsilon$-dependent. 
To prove our main result, i.e., Theorem \ref{thm1.2}, we need to compare $({\mathbf u}^{(0)},{\mathbf B^{(0)}},\theta^{(0)})$  with $({\mathbf u},{\mathbf B},\theta)$ and $({\mathbf u}^0,{\mathbf B^0},\theta^0)$, respectively.
It is clear that if the initial-layer corrections in (\ref{1.19})$_1$ are neglected, then the effective dynamics (\ref{1.19}) becomes the infinite Prandtl number dynamics (\ref{1.5}). It is natural to show that  as $\varepsilon\to0$, the solution  of the effective dynamics (\ref{1.19})--(\ref{1.20}) is close to the one  of the infinite Prandtl number dynamics  (\ref{1.5})--(\ref{1.6}).
However, the convergence from  $({\mathbf u},{\mathbf B},\theta)$ to $({\mathbf u}^{(0)},{\mathbf B^{(0)}},\theta^{(0)})$ is more complicated. Indeed, although $\tilde\uu|_{t=0}=(\uu-\uu^{(0)})|_{t=0}=0$, one cannot expect that  $\tilde\uu_t|_{t=0}=0$ also holds in general. As a result, the $H^2$-convergence of the velocity cannot be obtained directly. To circumvent this difficulty, we observe that the leading-order term of  temperature and Lorentz force (i.e., $A^{-1}\mathbb{P}(\kk \theta^{(0)}+\BB^{(0)}\cdot\nabla \BB^{(0)}+\BB^{(0)}_{y})$) plays an important role and its $t$-derivative acts as a correction term between $\uu_t$ and $\uu^{(0)}_t$ in some sense (see (\ref{3.40})--(\ref{3.43})).
The optimal convergence rates in (\ref{1.21}) and (\ref{1.22}) also need some careful analysis, based on the elementary  energy methods and the application of the Poincar${\rm\acute{e}}$'s inequality (see sections 3.2 and 3.3).

\section{Proof of Theorem \ref{thm1.1}}
The global existence of classical solutions $(\uu,\BB,\theta)$ to the problem \eqref{1.3}--\eqref{1.4} with smooth data can be easily proved via the standard Faedo-Galerkhin method and the global a priori estimates. The global solutions $(\uu^0,\BB^0,\theta^0)$ of  (\ref{1.5})--(\ref{1.6}) can be obtained as the vanishing $\varepsilon$-limit of $(\uu,\BB,\theta)$. Thus, for any given $0<T<\infty$, we assume that  $(\uu,\BB,\theta)$ and $(\uu^0,\BB^0,\theta^0)$ are   smooth solutions of \eqref{1.3}--\eqref{1.4} and (\ref{1.5})--(\ref{1.6}) on $\Omega \times[0, T]$, respectively. To prove Theorem \ref{thm1.1}, it suffices to derive some global (uniform-in-$\varepsilon$) estimates of $(\uu,\BB,\theta)$ and $(\uu^0,\BB^0,\theta^0)$.
\subsection{Global $\varepsilon$-independent estimates of $(\uu,\BB,\theta)$}
The purpose of this subsection is to derive the global uniform estimates of  $(\uu,\BB,\theta)$. For simplicity, throughout this paper we  use the same letter $C$  to denote the $\varepsilon$-independent constant.
\begin{pro}\label{pro2.1} Let $(\mathbf{u}, \mathbf{B},\theta)$ be a smooth solution of \eqref{1.3}--\eqref{1.4} on $\Omega \times[0, T]$. Then there exists a positive constant $C$, independent of $\varepsilon$, such that for any $p\geq 2$,
\begin{equation}\label{2.1}
\begin{split}
&\sup_{0\leq t\leq T} \|(\mathbf{B},\theta)(t)\|_{{H^1}{\cap} {W^{1,p}}}
+\int_{0}^{T}\left(\|(\mathbf{B},\theta)\|_{H^2}^2+\|(\mathbf{B}_t,\theta_t)\|_{L^2}^2\right)  dt\\
&\quad+\sup_{0\leq t\leq T}\|\mathbf{u}(t)\|_{H^2}
+\int_{0}^{T}\left(\|\mathbf{u}\|_{H^3}^2+\varepsilon\|\mathbf{u}_t\|_{H^1}^2\right)  dt \leq  C(p).
\end{split}
\end{equation}
\end{pro}
\pf The proofs are split into three steps.

{\bf Step I. the $L^2$-estimates}

First, it is easily deduced from (\ref{1.3})$_3$ and the maximum principle that
\begin{equation}\label{2.2}
\sup_{0\leq t \leq T}\|\theta(t)\|_{ L^\infty}\leq C.
\end{equation}

Let $\Theta\triangleq\theta-(1-y)$. Then, it holds that $\Theta|_{y=0,1}=0$ and
\begin{equation}\label{2.3}
\Theta_t+\uu\cdot\nabla\Theta=\Delta\Theta+u_2
\end{equation}
So, multiplying (\ref{1.3})$_1$, (\ref{1.3})$_2$ and (\ref{2.3}) by $\uu$, $\BB$ and $\Theta$ in $L^2$ respectively, integrating by parts, and using the Gronwall's inequality, we get that
\begin{equation}\label{2.4}
\sup_{0\leq t\leq T}\left(\varepsilon\|\uu\|_{L^2}^2+\|\BB\|_{L^2}^2+\|\Theta\|_{L^2}^2\right)(t)+\int_0^T\|(\nabla\uu,\nabla\BB,\nabla\Theta)\|_{L^2}^2dt\leq C.
\end{equation}

To prove the $\varepsilon$-independent estimate of $\|\uu\|_{L^2}$, we first multiply (\ref{1.3})$_2$ by $|\BB|^2\BB$ and integrate by parts to get
\begin{equation}\label{2.5}
\begin{split}
\frac{d}{dt}\|\BB\|_{L^4}^4+\||\BB||\nabla\BB|\|_{L^2}^2+\|\nabla|\BB|^2\|_{L^2}^2\leq& C\|\nabla \uu\|_{L^2}\left(\||\BB|^2\|_{L^4}^2+\|\BB\|_{L^4}\||\BB|^2\|_{L^4}\right)\\
\leq& C\left(1+\|\nabla \uu\|_{L^2}^2\right)\left(1+\|\BB\|_{L^4}^4\right)+\frac{1}{2}\|\nabla|\BB|^2\|_{L^2}^2,
\end{split}
\end{equation}
where we have used the Cauchy-Schwarz's inequality and the following Sobolev's inequality:
\begin{equation}\label{2.6}
\|f\|_{L^4}^2\leq C\|f\|_{L^2}^2+C\|f\|_{L^2}\|\nabla f\|_{L^2}
\end{equation}
for any $f\in \{f\in H^1 : f $ is $x$-periodic and $f=0$ on $ y=0,1 \}$.

Thus, using (\ref{2.4}) and the Gronwall's inequality, we infer from (\ref{2.5})  that
\begin{equation}\label{2.7}
\sup_{0\leq t\leq T}\|\BB(t)\|_{L^4}^4+\int_0^T\left(\||\BB||\nabla\BB|\|_{L^2}^2+\|\nabla|\BB|^2\|_{L^2}^2\right)dt\leq C.
\end{equation}

Now, multiplying (\ref{1.3})$_1$ by $\uu$ in $L^2$ again and integrating by parts, we find
\begin{equation*}
\varepsilon\frac{d}{dt}\|\uu\|_{L^2}^2+\|\nabla \uu\|_{L^2}^2\leq C\left(\|\uu\|_{L^2}\|\theta\|_{L^2}+\|\nabla \uu\|_{L^2}\|\BB\|_{L^4}^2+\|\nabla \uu\|_{L^2}\|\BB\|_{L^2}\right).
\end{equation*}
Thus, thanks to (\ref{2.4}), (\ref{2.7}) and the Poincar${\rm\acute{e}}$'s inequality:
$$
\|f\|_{L^2}\leq C\|\nabla f\|_{L^2}
$$
for any $f\in \{f\in H^1 : f $ is $x$-periodic and $f=0$ on $ y=0,1 \}$,
we have
\begin{equation*}
\varepsilon\frac{d}{dt}\|\uu\|_{L^2}^2+\|\uu\|_{L^2}^2\leq C,
\end{equation*}
and consequently,
$$
\frac{d}{dt}\left(e^{t/\varepsilon}\|\uu\|_{L^2}^2\right)\leq  C \varepsilon^{-1} e^{t/\varepsilon},
$$
which, integrated in time, shows that $\|\uu(t)\|_{L^2}$ is uniformly bounded in $\varepsilon$ for all $t\in[0,T]$.

To summarize, we have proved that
\begin{equation}\label{2.8}
\begin{split}
&\sup_{0\leq t\leq T} \left(\|(\uu,\BB,\theta)(t)\|_{L^2}^2+\|\BB(t)\|_{L^4}^4+\|\theta(t)\|_{L^\infty}\right)\\
&\qquad+\int_0^T\left(\|(\nabla\uu,\nabla\BB,\nabla\theta)\|_{L^2}^2+\||\BB||\nabla\BB|\|_{L^2}^2\right)dt\leq C.
\end{split}
\end{equation}

{\bf Step II. the $H^1$-estimates}

Note that
$$
\|(\Delta\uu,\Delta\BB,\Delta \Theta)\|_{L^2}=\|(\nabla^2\uu,\nabla^2\BB,\nabla^2\Theta)\|_{L^2}
$$
and
\begin{equation}\label{2.9}
\|\nabla f\|_{L^4}^2\leq C\|\nabla f\|_{L^2}^2+C\|\nabla f\|_{L^2}\|\nabla^2 f\|_{L^2},
\end{equation}
where  $f\in\{\uu,\BB,\Theta\}$.
Thus, multiplying \eqref{1.3}$_{1}$,  \eqref{1.3}$_{2}$ and \eqref{2.3} by $-\Delta \uu$, $-\Delta \BB$ and $-\Delta\Theta$ in $L^2$ respectively, and integrating by parts, we infer from (\ref{2.6}), (\ref{2.8}), (\ref{2.9}) and the Cauchy-Schwarz's inequality that
\begin{equation}\label{2.10}
\begin{split}
&\frac{d}{dt}\left(\varepsilon\|\nabla \uu\|_{L^2}^2+\|\nabla\BB\|_{L^2}^2+\|\nabla \Theta\|_{L^2}^2\right)+\|(\nabla^2\uu,\nabla^2\BB,\nabla^2\Theta)\|_{L^2}^2\\
&\quad\leq C\left(1+\|(\nabla\uu,\nabla\BB)\|_{L^2}^2\right)\left(1+\varepsilon\|\nabla \uu\|_{L^2}^2+\|\nabla\BB\|_{L^2}^2+\|\nabla \Theta\|_{L^2}^2\right),
\end{split}
\end{equation}
and hence, it follows from (\ref{2.8}), (\ref{2.10}) and the Gronwall's inequality that
\begin{equation}\label{2.11}
\sup_{0\leq t \leq T} \left(\varepsilon\|\nabla\uu\|_{L^2}^2+\|\nabla\BB\|_{L^2}^2+\|\nabla\Theta\|_{L^2}^2\right)(t)
+\int_{0}^{T}\|(\nabla^2\uu,\nabla^2\BB,\nabla^2\Theta)\|^2_{L^2}dt \leq C,
\end{equation}
which, together with (\ref{1.3}) and (\ref{2.3}), also gives
\begin{equation}\label{2.12}
\int_0^T\left(\varepsilon^2\|\uu_t\|_{L^2}^2+\|(\BB_t,\Theta_t)\|^2_{L^2}\right)dt\leq C.
\end{equation}

To prove the  $\varepsilon$-independent estimate of $\|\nabla\uu\|_{L^2}$, we multiply (\ref{1.3})$_1$ by $\uu_t$ in $L^2$ and integrate by parts to get that
\begin{equation}\label{2.13}
\begin{split}
\frac{d}{dt}\|\nabla\uu\|_{L^2}^2+\varepsilon\|\uu_t\|_{L^2}^2=&\langle\kk \theta+ \BB\cdot\nabla\BB+\BB_y,\uu_t\rangle-\varepsilon\langle\uu\cdot\nabla\uu,\uu_t \rangle\\
\triangleq& \frac{d}{dt}\langle\kk \theta+ \BB\cdot\nabla\BB+\BB_y, \uu\rangle+I,
\end{split}
\end{equation}
where $\langle\cdot,\cdot\rangle$ denotes the standard $L^2$-inner product and
$$
I\triangleq-\langle\kk \theta_t+ \BB_t\cdot\nabla\BB+\BB\cdot\nabla\BB_t+\BB_{ty},\uu\rangle-\varepsilon\langle\uu\cdot\nabla\uu,\uu_t \rangle.
$$

On one hand, by (\ref{2.8}) we have
\begin{equation*}
\langle\kk \theta+ \BB\cdot\nabla\BB+\BB_y, \uu\rangle=\langle\kk \theta,\uu\rangle-\langle \BB\cdot\nabla\uu,\BB\rangle-\langle \BB, \uu_y\rangle\leq\frac{1}{2}\|\nabla\uu\|_{L^2}^2+C.
\end{equation*}
On the other hand, based upon integration by parts, we have
\begin{equation*}
\begin{split}
I
&=-\langle\kk \theta_t,\uu\rangle+\langle \BB_t\cdot\nabla\uu,\BB\rangle+\langle\BB\cdot\nabla\uu,\BB_t\rangle+\langle\BB_{t},\uu_y\rangle-\varepsilon\langle\uu\cdot\nabla\uu,\uu_t \rangle\\
&\leq C\left(1+\|(\uu,\BB)\|_{H^2}^2\right)\|\nabla\uu\|_{L^2}^2+C\left(\|\theta_t\|_{L^2}^2+\|\BB_t\|_{L^2}^2+ \varepsilon^2\|\uu_t\|_{L^2}^2\right),
\end{split}
\end{equation*}
where have used  the Sobolev embedding inequality:
\begin{equation}\label{2.14}
\|f\|_{L^\infty}\leq C\|f\|_{H^2}\quad{\rm for\ \ any}\quad f\in H^2.
\end{equation}

Thus, inserting the above estimates into (\ref{2.13}), using (\ref{2.8}), (\ref{2.11}), (\ref{2.12}) and the Gronwall's inequality, we obtain
\begin{equation}\label{2.15}
\sup_{0\leq t\leq T}\|\nabla\uu(t)\|_{L^2}^2+\varepsilon\int_0^T\|\uu_t\|_{L^2}^2dt\leq C.
\end{equation}

{\bf Step III. the higher regularities}

As aforementioned, it is difficult to obtain the uniform $H^2$-estimate of $\BB$. Indeed, instead of this, we have the $W^{1,p}$-estimate of $\BB$ for any $p>2$, which particularly indicates that $\BB$ is uniformly bounded in $\varepsilon$. To do this, differentiating (\ref{1.3})$_2$ with respect to $y$, multiplying the resulting equation by $p|\BB_y|^{p-2}\BB_y$, and integrating by parts, we deduce that (noting that $(\BB_y+\uu)_y=0$ on $y=0,1$)
\begin{equation}
\begin{split}
&\frac{d}{dt}\|\BB_y\|_{L^p}^p+\int\left(|\BB_y|^{p-2}|\BB_{yy}|^2+|\BB_y|^{p-4}\left|\left(|\BB_y|^2\right)_y\right|^2\right) dxdy\\
&\qquad+\int\left(|\BB_y|^{p-2}|\BB_{xy}|^2+|\BB_y|^{p-4}\left|\left(|\BB_y|^2\right)_x\right|^2\right)dxdy\\
&\quad\leq  C(p)\int \left(|\uu|^2|\nabla\BB|^2+|\BB|^2|\nabla\uu|^2+|\uu_y|^2\right)|\BB_y|^{p-2}dxdy\\
&\quad\leq C(p)\left(\|(\uu,\BB)\|_{H^1}^2\|(\nabla\uu,\nabla\BB)\|_{H^1}^2+\|\nabla\uu\|_{H^1}^2\right)\|\BB_y\|_{L^p}^{p-2}\\
&\quad\leq  C(p)\left(\|\nabla\BB\|_{H^1}^2+\|\nabla\uu\|_{H^1}^2\right)\left(1+\|\BB_y\|_{L^p}^{p}\right),
\end{split}\label{2.16}
\end{equation}
where we have used (\ref{2.11}), (\ref{2.15}) and the Sobolev embedding inequality $H^1\hookrightarrow L^q$ for any $q\geq1$. So, using (\ref{2.11}) and the Gronwall's inequality, we have from  (\ref{2.16}) that
\begin{equation}
\sup_{0\leq t\leq T}\|\BB_y(t)\|_{L^p}\leq C(p)\quad{\rm for}\quad \forall \ p\geq 2,\label{2.17}
\end{equation}
and similarly,
\begin{equation*}
\sup_{0\leq t\leq T}\|\BB_x(t)\|_{L^p}\leq C(p)\quad{\rm for}\quad \forall \ p\geq 2,
\end{equation*}
so that, it follows from the Sobolev's embedding inequality $W^{1,p}\hookrightarrow L^\infty$ for $p>2$  that
\begin{equation}
\|\BB(t)\|_{L^\infty}+\|\BB(t)\|_{W^{1,p}}\leq C(p)\quad{\rm for}\quad \forall\ t\in[0,T].\label{2.18}
\end{equation}

Analogously to the proof of (\ref{2.17}), one also gets that for any $t\in[0,T]$,
\begin{equation}
\|\Theta(t)\|_{W^{1,p}}+\|\theta(t)\|_{W^{1,p}}\leq C(p)\quad{\rm for}\quad \forall\ p>2.\label{2.19}
\end{equation}

The estimate of $\|\uu\|_{H^2}$ needs more works. For this purpose, we first
differentiate (\ref{1.3})$_1$ with respect to $t$, multiply the resulting equation by $\uu_t$ in $L^2$, and integrate by parts to deduce that
\begin{equation*}
\begin{split}
\frac{\varepsilon}{2}\frac{d}{dt}\|\uu_t\|_{L^2}^2+\|\nabla\uu_t\|_{L^2}^2=&-\varepsilon\langle \uu_t\cdot\nabla\uu,\uu_t\rangle+\langle\kk\theta_t,\uu_t\rangle+\langle\BB_{yt},\uu_{t}\rangle
\\
&+\langle\BB\cdot\nabla\BB_t,\uu_t\rangle+\langle\BB_t\cdot\nabla\BB,\uu_t\rangle\triangleq II.
\end{split}
\end{equation*}
Using (\ref{2.14}), (\ref{2.18}) and the Poincar${\rm\acute e}$'s inequality, we obtain after integrating by parts that
\begin{equation*}
\begin{split}
II=&\varepsilon\langle \uu_t\cdot\nabla\uu_t,\uu\rangle+\langle\kk\theta_t,\uu_t\rangle-\langle\BB_t,\uu_{yt}\rangle
-\langle\BB\cdot\nabla\uu_t,\BB_t\rangle-\langle\BB_t\cdot\nabla\uu_t,\BB\rangle\\
\leq& \frac{1}{2}\|\nabla\uu_t\|_{L^2}^2+C \varepsilon\|\uu_t\|_{L^2}^2\|\uu\|_{H^2}^2+C \|(\BB_t,\theta_t)\|_{L^2}^2,
\end{split}
\end{equation*}
so that
\begin{equation}\label{2.20}
 \varepsilon \frac{d}{dt}\|\uu_t\|_{L^2}^2+\|\nabla\uu_t\|_{L^2}^2\leq C \varepsilon\|\uu_t\|_{L^2}^2\|\uu\|_{H^2}^2+C \|(\BB_t,\theta_t)\|_{L^2}^2.
\end{equation}

In a similar manner, we also have
\begin{equation*}
\begin{split}
\frac{1}{2}\frac{d}{dt}\|(\BB_t,\theta_t)\|_{L^2}^2+\|(\nabla\BB_t,\nabla\theta_t)\|_{L^2}^2=&\langle\BB\cdot\nabla\uu_t,\BB_t\rangle +\langle\BB_t\cdot\nabla\uu,\BB_t\rangle+\langle\uu_{yt},\BB_t\rangle\\
&-\langle\uu_t\cdot\nabla\BB,\BB_t\rangle-\langle\uu_t\cdot\nabla\theta,\theta_t\rangle
\triangleq III,
\end{split}
\end{equation*}
where we can utilize (\ref{2.2}), (\ref{2.14}), (\ref{2.18}) and the Poincar${\rm\acute{e}}$ inequality to get that
\begin{equation*}
\begin{split}
III=&\langle\BB\cdot\nabla\uu_t,\BB_t\rangle -\langle\BB_t\cdot\nabla\BB_t,\uu\rangle+\langle\uu_{yt},\BB_t\rangle+\langle\uu_t\cdot\nabla\BB_t,\BB\rangle+\langle\uu_t\cdot\nabla\theta_t,\theta\rangle\\
\leq& \frac{1}{2}\|(\nabla\BB_t,\nabla\theta_t)\|_{L^2}^2+ C\left(1+\|\uu\|_{H^2}^2\right) \|(\BB_t,\theta_t)\|_{L^2}^2+C\|\nabla \uu_t\|_{L^2}^2,
\end{split}
\end{equation*}
and consequently,
\begin{equation}\label{2.21}
\frac{d}{dt}\|(\BB_t,\theta_t)\|_{L^2}^2+\|(\nabla\BB_t,\nabla\theta_t)\|_{L^2}^2\leq C\left(1+\|\uu\|_{H^2}^2\right) \|(\BB_t,\theta_t)\|_{L^2}^2+C\|\nabla \uu_t\|_{L^2}^2.
\end{equation}

As a result of (\ref{2.20}) and (\ref{2.21}), we find
\begin{equation}\label{2.22}
\begin{split}
&\frac{d}{dt}\left(\varepsilon\|\uu_t\|_{L^2}^2+\|(\BB_t,\theta_t)\|_{L^2}^2\right)+\|(\nabla\uu_t,\nabla\BB_t,\nabla\theta_t)\|_{L^2}^2\\
&\quad \leq C\left(1+\|\uu\|_{H^2}^2\right)\left(\varepsilon\|\uu_t\|_{L^2}^2+\|(\BB_t,\theta_t)\|_{L^2}^2\right).
\end{split}
\end{equation}

To eliminate the effect of initial layer, multiplying (\ref{2.22}) by $t$ and integrating it over $(0,T)$, we deduce from  (\ref{2.12}), (\ref{2.15}) and the Gronwall's inequality that
\begin{equation}\label{2.23}
\begin{split}
&\sup_{0\leq t\leq T}\left[t\left(\varepsilon\|\uu_t(t)\|_{L^2}^2+\|(\BB_t,\theta_t)(t)\|_{L^2}^2\right)\right]+\int_0^Tt \|(\nabla\uu_t,\nabla\BB_t,\nabla\theta_t)\|_{L^2}^2dt\\
&\quad \leq C\int_0^T\left(\varepsilon\|\uu_t\|_{L^2}^2+\|(\BB_t,\theta_t)\|_{L^2}^2\right)dt\leq C.
\end{split}
\end{equation}

Since $\uu_0\in H^2$ implies that $\varepsilon\uu_t|_{t=0}\in L^2$, after multiplying (\ref{2.22}) by $\varepsilon$ and integrating it over $(0,T)$, we arrive at
\begin{equation}\label{2.24}
\sup_{0\leq t\leq T}\left(\varepsilon^2 \|\uu_t\|_{L^2}^2+\varepsilon\|(\BB_t,\theta_t)\|_{L^2}^2\right)(t) +\varepsilon\int_0^T\|(\nabla\uu_t,\nabla\BB_t,\nabla\theta_t) \|_{L^2}^2dt\leq C.
\end{equation}

In view of (\ref{2.11}), (\ref{2.18}) and (\ref{2.24}), we can make use of   (\ref{2.6}), (\ref{2.9}) and the standard estimates of the Stokes equations to get from (\ref{1.3})$_1$ that
\begin{equation*}
\begin{split}
\|\nabla^2\uu\|_{L^2}\leq& C\varepsilon\left(\|\uu_t\|_{L^2}+\|\uu\|_{L^4}\|\nabla\uu\|_{L^4}\right)+C\left(\|\theta\|_{L^2}+\|\BB\|_{L^\infty}\|\nabla\BB\|_{L^2}+\|\nabla\BB\|_{L^2}\right)\\
\leq& C+C\|\nabla^2\uu\|_{L^2}^{1/2}\leq \frac{1}{2}\|\nabla^2\uu\|_{L^2}+C,
\end{split}
\end{equation*}
which immediately results in
\begin{equation}
\sup_{0\leq t\leq T}\|\nabla^2\uu(t)\|_{L^2}\leq C,\label{2.25}
\end{equation}
and moreover,
\begin{equation}\label{2.26}
\begin{split}
\int_0^T\|\uu\|_{H^3}^2dt\leq& C+C\int_0^T\left(\varepsilon^2\|\uu_t\|_{H^1}^2+\|(\uu,\BB)\|_{L^\infty}^2\|(\uu,\BB)\|_{H^2}^2\right)dt\\
&+C\int_0^T\left(\|(\nabla\uu,\nabla\BB)\|_{L^4}^4+\|(\BB,\theta)\|_{H^2}^2\right)dt\\
\leq& C.
\end{split}
\end{equation}

Now, collecting (\ref{2.8}), (\ref{2.11}), (\ref{2.12}), (\ref{2.15}), (\ref{2.18}), (\ref{2.19}) and (\ref{2.23})--(\ref{2.26}) together leads to the desired estimates stated in  (\ref{2.1}). The proof of Proposition \ref{pro2.1} is therefore complete.\hfill$\square$

\subsection{Global estimates of $(\uu^0,\BB^0,\theta^0)$}
 This subsection is devoted to  the global estimates of the solutions to the problem \eqref{1.5}--\eqref{1.6}, which can be achieved via the standard estimates of the Stokes equations.
\begin{pro}\label{pro2.2} Let $(\mathbf{u}^0, \mathbf{B}^0,\theta^0)$ be a smooth solution of \eqref{1.5}--\eqref{1.6} on $\Omega \times[0, T]$. Then there exists a positive constant $C$, such that
\begin{equation}\label{2.27}
\sup_{0\leq t\leq T}\left(\|(\mathbf{u}^0,\mathbf{B}^0,\theta^0)\|_{H^2}^2+\|(\mathbf{B}_t^0,\theta_t^0)\|_{L^2}^2\right)(t)
+\int_{0}^{T}\|(\mathbf{u}^0_t,\mathbf{B}_t^0,\theta_t^0)\|_{H^1}^2  dt \leq  C.
\end{equation}
\end{pro}
\pf First, one easily gets from  (\ref{1.5})  that
\begin{equation}\label{2.28}
\sup_{0\leq t\leq T}\left(\|(\BB^0,\theta^0)(t)\|_{L^2}^2+\|\theta^0(t)\|_{L^\infty}\right)
+\int_{0}^{T}\|(\nabla\uu^0,\nabla\BB^0,\nabla\theta^0)\|_{L^2}^2  dt \leq  C.
\end{equation}

Next, similarly to the proof of (\ref{2.7}), by  (\ref{2.28}) we have
\begin{equation}\label{2.29}
\sup_{0\leq t\leq T}\|\BB^0(t)\|_{L^4}^4+\int_0^T\left(\||\BB^0||\nabla\BB^0|\|_{L^2}^2+\|\nabla|\BB^0|^2\|_{L^2}^2\right)dt\leq C.
\end{equation}
Thus, it is easily seen from (\ref{1.5})$_1$ and (\ref{2.29}) that
\begin{equation}\label{2.30}
\|\nabla \uu^0(t)\|_{L^2}\leq C\left(\|\BB^0(t)\|_{L^4}^2+\|(\BB^0,\theta^0)(t)\|_{L^2}\right)\leq C,\quad \forall\ t\in[0,T],
\end{equation}
and moreover, it follows from the $H^2$-estimate of the Stokes equations that
\begin{equation}\label{2.31}
\int_0^T\|\nabla^2\uu^0\|_{L^2}^2dt\leq C\int_0^T\left(\||\BB^0||\nabla\BB^0|\|_{L^2}^2+\|\nabla \BB^0\|_{L^2}^2+\|\theta^0\|_{L^2}^2\right)dt\leq C.
\end{equation}

With the help of (\ref{2.28})--(\ref{2.31}), we can obtain in a manner similar to the derivations of (\ref{2.11}) and (\ref{2.12}) that
\begin{equation}\label{2.32}
\sup_{0\leq t\leq T} \|(\nabla\BB^0,\nabla\theta^0)(t)\|_{L^2}^2+\int_0^T\left(\|(\nabla^2\BB^0,\nabla^2\theta^0)\|_{L^2}^2+\|(\BB^0_t,\theta^0_t)\|_{L^2}^2\right)dt\leq C.
\end{equation}

Finally, differentiating (\ref{1.5})$_1$, (\ref{1.5})$_2$ and (\ref{1.5})$_3$ with respect to $t$ and multiplying the resulting equations by $\uu_t$, $\BB_t$ and $\theta_t$ in $L^2$ respectively, and integrating by parts, we have
\begin{equation}\label{2.33}
\begin{split}
&\frac{d}{dt}\left(\|\BB_t^0\|_{L^2}^2+\|\theta_t^0\|_{L^2}^2\right)+\|(\nabla\uu^0_t,\nabla\BB^0_t,\nabla\theta^0_t)\|_{L^2}^2\\
&\quad \leq C\left(1+\|(\uu^0,\BB^0,\theta^0)\|_{H^2}^2\right)\left( \|\BB^0_t\|_{L^2}^2+\|\theta^0_t\|_{L^2}^2\right),
\end{split}
\end{equation}
where we have used (\ref{2.6}), (\ref{2.28})--(\ref{2.32})  and the Poincar${\rm\acute{e}}$'s inequality to get that
\begin{equation*}
\begin{split}
\langle\kk\theta^0_t,\uu^0_t\rangle+\langle\partial_y\BB^0_t,\uu^0_t\rangle +\langle\partial_y\uu^0_t,\BB^0_t\rangle
&\leq \|\uu^0_t\|_{L^2}\|\theta^0_t\|_{L^2}+ \|\nabla\uu^0_t\|_{L^2}\|\BB^0_t\|_{L^2}\\
&\leq\frac{1}{8}\|\nabla\uu^0_t\|_{L^2}^2+C\left(\|\BB^0_t\|_{L^2}^2+\|\theta^0_t\|_{L^2}^2\right),
\end{split}
\end{equation*}
\begin{equation*}
\begin{split}
&-\langle(\uu^0\cdot\nabla\BB^0)_t,\BB^0_t\rangle-\langle(\uu^0\cdot\nabla\theta^0)_t,\theta^0_t\rangle =-\langle \uu^0_t\cdot\nabla\BB^0,\BB^0_t\rangle-\langle \uu^0_t\cdot\nabla\theta^0,\theta^0_t\rangle\\
&\quad\leq \|\uu^0_t\|_{L^4}\left(\|\nabla\BB^0\|_{L^4}+\|\nabla\theta^0\|_{L^4}\right)\left(\|\BB^0_t\|_{L^2}+\|\theta^0_t\|_{L^2}\right)\\
&\quad\leq \frac{1}{8} \|\nabla\uu^0_t\|_{L^2} ^2+C\left(\|\nabla\BB^0\|_{H^1}^2+\|\nabla\theta^0\|_{H^1}^2\right) \left(\|\BB^0_t\|_{L^2}^2+\|\theta^0_t\|_{L^2}^2\right),
\end{split}
\end{equation*}
and
\begin{equation*}
\begin{split}
&\langle(\BB^0\cdot\nabla\BB^0)_t,\uu^0_t\rangle+\langle(\BB^0\cdot\nabla\uu^0)_t,\BB^0_t\rangle=\langle\BB^0_t\cdot\nabla\BB^0,\uu^0_t\rangle+\langle\BB^0_t\cdot\nabla\uu^0,\BB^0_t\rangle\\
&\quad \leq \|\BB^0_t\|_{L^4}\|\nabla\BB^0\|_{L^4}\|\uu^0_t\|_{L^2}+ \|\BB^0_t\|_{L^4}\|\nabla\uu^0\|_{L^4}\|\BB^0_t\|_{L^2} \\
&\quad\leq C\left(\|\BB^0_t\|_{L^2}+\|\BB^0_t\|_{L^2}^{1/2}\|\nabla \BB^0_t\|_{L^2}^{1/2}\right)
\left(\|\nabla\BB^0\|_{L^2}+\|\nabla\BB^0\|_{L^2}^{1/2}\|\nabla^2\BB^0\|_{L^2}^{1/2}\right)\| \uu^0_t\|_{L^2}\\
&\qquad +C\left(\|\BB^0_t\|_{L^2}+\|\BB^0_t\|_{L^2}^{1/2}\|\nabla \BB^0_t\|_{L^2}^{1/2}\right)
\left(\|\nabla\uu^0\|_{L^2}+\|\nabla\uu^0\|_{L^2}^{1/2}\|\nabla^2\uu^0\|_{L^2}^{1/2}\right)\|\BB^0_t\|_{L^2}\\
&\quad \leq \frac{1}{8}\|(\nabla\uu^0_t,\nabla\BB^0_t)\|_{L^2}^2+C\left(1+\|(\nabla^2\uu^0,\nabla^2\BB^0)\|_{L^2}^2\right)\|\BB^0_t\|_{L^2}^2.
\end{split}
\end{equation*}

Using (\ref{2.28}), (\ref{2.31}), (\ref{2.32}) and the Gronwall's inequality, we deduce from (\ref{2.33})  that
\begin{equation}\label{2.34}
\sup_{0\leq t\leq T}\left(\|\BB_t^0\|_{L^2}^2+\|\theta_t^0\|_{L^2}^2\right)(t)+\int_0^T\|(\nabla\uu^0_t, \nabla\BB^0_t, \nabla\theta^0_t)\|_{L^2}^2dt\leq C.
\end{equation}
Thus, it follows from (\ref{1.5})$_{2,3}$, (\ref{2.32}), (\ref{2.34}) and the Sobolev embedding inequality that
\begin{equation}\label{2.35}
\begin{split}
\|(\BB^0,\theta^0)(t)\|_{H^2}\leq& C+C\left(\|\BB^0_t\|_{L^2}+\|\theta^0_t\|_{L^2}+\|\uu^0\|_{H^1}\right)\\
&+C\left(\||\uu^0||\nabla\BB^0|\|_{L^2}+\||\BB^0||\nabla\uu^0|\|_{L^2}+\||\uu^0||\nabla\theta^0|\|_{L^2}\right)\\
\leq&C+C\|\uu^0\|_{H^2},
\end{split}
\end{equation}
and similarly, it follows from (\ref{1.5})$_1$ and the standard estimates of the Stokes equations that
\begin{equation}\label{2.36}
\begin{split}
\|\uu^0(t)\|_{H^2}\leq& C+C\left(\||\BB^0||\nabla\BB^0|\|_{L^2}+\|\BB^0\|_{H^1}+\|\theta^0\|_{L^2}\right)\\
\leq& C+C\left(\|\nabla\BB^0\|_{L^2}+\|\nabla\BB^0\|_{L^2}^{1/2}\|\nabla^2\BB^0\|_{L^2}^{1/2}\right)\\
\leq& C+C\|\BB^0\|_{H^2}^{1/2}.
\end{split}
\end{equation}

As a result, we conclude from   (\ref{2.35}), (\ref{2.36}) and  the Cauchy-Schwarz's inequality that
$$
\sup_{0\leq t\leq T}\|(\uu^0,\BB^0,\theta^0)(t)\|_{H^2}\leq C.
$$
This, together with (\ref{2.28}), (\ref{2.32}) and (\ref{2.34}), finishes the proof of Proposition \ref{pro2.2}.\hfill$\square$

\subsection{Convergence from $(\uu,\BB,\theta)$ to $(\uu^0,\BB^0,\theta^0)$}
With the help of the global (uniform) estimates stated in Propositions \ref{pro2.1} and \ref{pro2.2}, we are now ready to prove  Theorem \ref{thm1.1}. First of all, the global existence of strong solutions to the problem (\ref{1.3})--(\ref{1.4}) is an immediate consequence of the global estimates in Proposition \ref{pro2.1}, and the global solutions of  (\ref{1.5})--(\ref{1.6}) can be obtained as the vanishing $\varepsilon$-limit of $(\uu,\BB,\theta)$. So, it only remains to prove the convergence rates. To do this, we define
$$
\bar\uu=\uu-\uu^0,\quad\bar\BB=\BB-\BB^0,\quad\bar\theta=\theta-\theta^0.
$$
Then it is easily derived from (\ref{1.3})--(\ref{1.4}) and (\ref{1.5})--(\ref{1.6}) that
\begin{equation}\label{2.37}
\begin{split}
\varepsilon\left(\frac{\partial \bar\uu}{\partial t}+\uu\cdot \nabla\bar\uu\right)+\nabla q-\Delta\bar\uu=&\kk \bar\theta+ \bar\BB\cdot \nabla \BB+\BB^0\cdot \nabla \bar\BB+\frac{\partial\bar\BB}{\partial y}\\
&-\varepsilon\left(\uu^0_t+\uu^0\cdot \nabla\uu^0+\bar\uu\cdot \nabla \uu^0\right),
\end{split}
\end{equation}
\begin{equation}\label{2.38}
\frac{\partial \bar\BB}{\partial t}+\uu\cdot \nabla \bar\BB-\Delta\bar\BB=-\bar\uu\cdot \nabla \BB^0+\bar\BB\cdot\nabla \uu+\BB^0\cdot \nabla\bar\uu+\frac{\partial \bar\uu}{\partial y},
\end{equation}
and
\begin{equation}\label{2.39}
\frac{\partial \bar\theta}{\partial t}+\uu\cdot \nabla \bar\theta-\Delta\bar\theta=-\bar\uu\cdot \nabla \theta^0,
\end{equation}
with $\divg\bar\uu=\divg\bar\BB=0$ and the vanishing initial-boundary conditions:
$$
 (\bar\BB,\bar\theta)|_{t=0}=0\quad {\rm and}\quad (\bar\uu,\bar\BB,\bar\theta)|_{y=0,1}=0.
$$

First, multiplying (\ref{2.37}) by $\bar\uu$ in $L^2$ and integrating by parts, we have from Propositions \ref{pro2.1} and \ref{pro2.2} that
\begin{equation*}
\begin{split}
\varepsilon\frac{d}{dt} \|\bar\uu\|_{L^2}^2+ \|\nabla\bar\uu\|_{L^2}^2
\leq& C\left(\|\bar\theta\|_{L^2}^2+\|(\BB,\BB^0)\|_{L^\infty}^2\|\bar\BB\|_{L^2}^2+\|\bar\BB\|_{L^2}^2\right)\\
&+C\varepsilon\left(\|\uu^0_t\|_{L^2}^2+\|\uu^0\|_{H^2}^4+\|\uu^0\|_{L^\infty}^2\|\bar\uu\|_{L^2}^2\right)\\
\leq&C\varepsilon+C\left(\varepsilon\|\bar\uu\|_{L^2}^2+\|\bar\BB\|_{L^2}^2+\|\bar\theta\|_{L^2}^2\right),
\end{split}
\end{equation*}
and analogously,
\begin{equation*}
\frac{d}{dt}\left(\|\bar\BB\|_{L^2}^2 + \|\bar\theta\|_{L^2}^2\right) +\|\nabla\bar\BB\|_{L^2}^2 +\|\nabla\bar\theta\|_{L^2}^2
\leq C\left(\|\bar\uu\|_{L^2}^2+\|\bar\BB\|_{L^2}^2+\|\bar\theta\|_{L^2}^2\right).
\end{equation*}

Owing to the Poincar${\rm\acute{e}}$'s inequality, it holds that $\|\tilde\uu\|_{L^2}\leq C\|\nabla\tilde\uu\|_{L^2}$. Thus,
\begin{equation*}
\begin{split}
&\frac{d}{dt} \left(\varepsilon\|\bar\uu\|_{L^2}^2+\|\bar\BB\|_{L^2}^2 + \|\bar\theta\|_{L^2}^2\right) + \|\nabla\bar\uu\|_{L^2}^2+\|\nabla\bar\BB\|_{L^2}^2 +\|\nabla\bar\theta\|_{L^2}^2\\
&\quad
\leq C \left(\varepsilon\|\bar\uu\|_{L^2}^2+\|\bar\BB\|_{L^2}^2+\|\bar\theta\|_{L^2}^2\right)
+C\varepsilon ,
\end{split}
\end{equation*}
so that
\begin{equation}
\sup_{0\leq t\leq T}\left(\varepsilon\|\bar\uu\|_{L^2}^2+\|\bar\BB\|_{L^2}^2 + \|\bar\theta\|_{L^2}^2\right)(t) + \int_0^T \|(\nabla\bar\uu,\nabla\bar\BB,\nabla\bar\theta)\|_{L^2}^2dt\leq C\varepsilon.\label{2.40}
\end{equation}

Multiplying (\ref{2.38}) and (\ref{2.39}) by $\bar\BB_t$ and $\bar\theta_t$ in $L^2$ respectively, integrating by parts, using Propositions \ref{pro2.1}--\ref{pro2.2} and the Poincar${\rm\acute{e}}$'s inequality, we  deduce
\begin{equation*}
\begin{split}
&\frac{d}{dt}\left(\|\nabla\bar\BB\|_{L^2}^2+\|\nabla\bar\theta\|_{L^2}^2\right)+\|\bar\BB_t\|_{L^2}^2+\|\bar\theta_t\|_{L^2}^2\\
&\quad\leq C\|\uu\|_{H^2}^2\left(\|\nabla\bar\BB\|_{L^2}^2+\|\nabla\bar\theta\|_{L^2}^2\right)+\left(1+\|\BB^0\|_{H^2}^2+\|\theta^0\|_{H^2}^2\right)\|\nabla\bar\uu\|_{L^2}^2\\
&\quad\leq C \left(\|\nabla\bar\BB\|_{L^2}^2+\|\nabla\bar\theta\|_{L^2}^2+\|\nabla\bar\uu\|_{L^2}^2\right),
\end{split}
\end{equation*}
and hence, by (\ref{2.40}) we have
\begin{equation}
\sup_{0\leq t\leq T}
\left(\|\nabla\bar\BB\|_{L^2}^2+\|\nabla\bar\theta\|_{L^2}^2\right)(t)+\int_0^T\left(\|\bar\BB_t\|_{L^2}^2+\|\bar\theta_t\|_{L^2}^2\right)dt\leq C\varepsilon,\label{2.41}
\end{equation}
which, together with (\ref{2.38}) and (\ref{2.39}), also yields
\begin{equation}\label{2.42}
\int_0^T\left(\|\nabla^2\bar\BB\|_{L^2}^2+\|\nabla^2\bar\theta\|_{L^2}^2\right)dt\leq C\varepsilon.
\end{equation}

Applying $\nabla$ to both sides of (\ref{2.38}), multiplying it by $|\nabla\bar\BB|^2\nabla\bar\BB$, and integrating by parts, we deduce in a manner similar to the derivation of (\ref{2.17}) that
\begin{equation}
\begin{split}
&\frac{d}{dt}\|\nabla\bar\BB\|_{L^4}^4+\||\nabla\bar\BB||\nabla^2\bar\BB|\|_{L^2}^2\\
&\quad\leq C\int \left(|\uu|^2|\nabla \bar\BB|^2+|\bar\uu|^2| \nabla \BB^0|^2+|\bar\BB|^2|\nabla \uu|^2+|\BB^0|^2|\nabla\bar\uu|^2+|\bar\uu_y|^2\right)|\nabla\bar \BB|^2dxdy\\
&\quad\leq C\left(\|\uu\|_{L^\infty}^2\|\nabla\bar\BB\|_{L^4}^2+\|\bar\uu\|_{H^1}^2\|\BB^0\|_{H^2}^2
+\|\bar\BB\|_{H^1}^2\|\nabla\uu\|_{H^1}^2 +\|\bar\uu_y\|_{L^4}^2\right)\|\nabla\bar\BB\|_{L^4}^2\\
&\quad\leq C\|\nabla\bar\BB\|_{L^4}^4+C\left(\varepsilon+\|\nabla\bar\uu\|_{L^2}^2\right),
\end{split}\label{2.43}
\end{equation}
where we have used (\ref{2.41}), Propositions \ref{pro2.1}--\ref{pro2.2}, the Sobolev embedding inequality, the Poincar${\rm\acute{e}}$'s inequality, and the following simple fact (noting that $\bar\uu=\uu-\uu^0\in L^\infty(0,T;H^2)$)
$$
\|\bar\uu_y\|_{L^4}^4\leq \|\bar\uu_y\|_{L^2}^4+\|\bar\uu_y\|_{L^2}^2\|\nabla\bar\uu_y\|_{L^2}^2\leq C\|\nabla\bar \uu\|_{L^2}^2.
$$

Thus, it follows from (\ref{2.40}), (\ref{2.43}) and the Gronwall's inequality that
\begin{equation}\label{2.44}
\sup_{0\leq t\leq T}\|\nabla\bar\BB(t)\|_{L^4}^4\leq C\varepsilon.
\end{equation}
In the exactly same way, we also have
\begin{equation}
\sup_{0\leq t\leq T}\|\nabla\bar\theta(t)\|_{L^4}^4\leq C\varepsilon.\label{2.45}
\end{equation}

Combining (\ref{2.40}), (\ref{2.41}), (\ref{2.44}), (\ref{2.45}) and the Sobolev embedding inequality, we see that
\begin{equation}\label{2.46}
(\BB,\theta)\to(\bar\BB,\bar\theta)\quad{\rm in}\quad C(\overline\Omega\times[0,T]),
\end{equation}
which particularly implies that there is no initial-layer between $(\BB,\theta)$ and $(\bar\BB,\bar\theta)$ in the sense of uniform convergence.

With the help of (\ref{2.41}) and Propositions \ref{pro2.1}--\ref{pro2.2}, it is easy to derive the convergence of $\bar\uu$ strictly away from the initial layer. Indeed, if rewriting (\ref{2.37}) in the form:
$$
\nabla q-\Delta\bar\uu=-\varepsilon(\uu_t+\uu\cdot\nabla\uu)+\kk \bar\theta+ \BB\cdot \nabla \bar\BB+\bar\BB\cdot \nabla \BB^0+\bar\BB_y,
$$
then we can utilize the estimates of the Stokes equations to infer from (\ref{2.40}), (\ref{2.41}) and Proposition \ref{pro2.1}--\ref{pro2.2} that
\begin{equation*}
\|\bar\uu\|_{H^2}^2\leq \varepsilon^2\left(\|\uu_t\|_{L^2}^2+\|\uu\|_{H^2}^4\right)+\|(\bar\BB,\bar\theta)\|_{H^1}^2
\leq C\varepsilon+C\varepsilon^2\|\uu_t\|_{L^2}^2,
\end{equation*}
which, combined with (\ref{2.23}), yields
\begin{equation}\label{2.47}
t\|\bar\uu(t)\|_{H^2}^2
\leq C\varepsilon+C\varepsilon^2 t\|\uu_t(t)\|_{L^2}^2\leq C\varepsilon\quad {\rm for}\quad\forall\ t\in[0,T].
\end{equation}

As a result, it follows from (\ref{2.47}) and the Sobolev embedding inequality that
\begin{equation}\label{2.48}
\|\bar\uu(t)\|_{C({\overline\Omega})}^2\leq C\|\bar\uu(t)\|_{H^2}^2\leq \frac{C\varepsilon}{t}\to 0,\quad {\rm if}\quad t\geq \varepsilon^{1-\alpha}\quad{\rm with}\quad\forall\ \alpha\in(0,1),
\end{equation}
which indicates that as $\varepsilon\to0$, $\uu$ converges to $\uu^0$ in $H^2$ strictly away from the initial layer, whose width is of the value $O(\varepsilon^{1-\alpha})$ with any $\alpha\in(0,1)$.

Collecting (\ref{2.40})--(\ref{2.42}) and (\ref{2.44})--(\ref{2.48}) together leads to the convergence results stated in (\ref{1.9}) and (\ref{1.10}). The proof of Theorem \ref{thm1.1} is therefore complete.\hfill$\square$

\section{Proof of Theorem \ref{thm1.2}}
In this section, we aim to prove Theorem \ref{thm1.2} by comparing the solutions $(\uu^{(0)},\BB^{(0)},\theta^{(0)})$ of the problem (\ref{1.19})--(\ref{1.20}) with the ones of the problems (\ref{1.5})--(\ref{1.6}) and (\ref{1.3})--(\ref{1.4}) successively.
\subsection{Global estimates of $(\uu^{(0)},\BB^{(0)},\theta^{(0)})$}
This subsection is devoted to the global estimates of $(\uu^{(0)},\BB^{(0)},\theta^{(0)})$. To this end, we first recall some known facts of the Stokes operator $A$. It is well known that (cf. \cite{L1963,T2001})
the operator $A$ is coercive and satisfies
\begin{equation}
 \|\uu\|_{L^2}\leq C\|\nabla\uu\|_{L^2}\leq C\|A\uu\|_{L^2},\label{3.1}
\end{equation}
and the semigroup $e^{-tA}$ satisfies
\begin{equation}
\|e^{-tA}\uu\|_{L^2}\leq Ce^{-t}\|\uu\|_{L^2}.\label{3.2}
\end{equation}
Moreover, if $\uu=A^{-1}\ff$  is a solution of the problem (\ref{1.14}), then it follows from (\ref{3.1}) that
\begin{equation}\label{3.3}
\|\nabla\uu\|_{L^2}^2=\langle A\uu,\uu\rangle=\langle \ff,\uu\rangle\leq \|\ff\|_{H^{-1}}\|\uu\|_{H^1}\leq C\|\ff\|_{H^{-1}}\|\nabla\uu\|_{L^2}
\end{equation}
and consequently,
\begin{equation}
\|A^{-1}\ff\|_{L^2}=\|\uu\|_{L^2}\leq C\|\nabla\uu\|_{L^2}\leq C\|\ff\|_{H^{-1}}.\label{3.4}
\end{equation}

With the help of (\ref{3.1})--(\ref{3.4}) and the estimates of Stokes equations, we can prove  that
\begin{pro}
\label{pro3.1}Let $(\mathbf{u}^{(0)}, \mathbf{B}^{(0)},\theta^{(0)})$ be a smooth solution of \eqref{1.19}--\eqref{1.20} on $\Omega \times[0, T]$. Then there exists a positive constant $C$, such that for any $p>2$,
\begin{equation}\label{3.5}
\begin{split}
&\sup_{0\leq t\leq T}\left(\|\mathbf{u}^{(0)}(t)\|_{H^2}+ \|(\mathbf{B}^{(0)},\theta^{(0)})(t)\|_{{H^1}{\cap} {W^{1,p}}}\right)\\
&\qquad
+\int_{0}^{T}\left(\|(\mathbf{B}^{(0)},\theta^{(0)})(t)\|_{H^2}^2+\|(\mathbf{B}^{(0)}_t,\theta_t^{(0)})\|_{L^2}^2\right)  dt \leq  C(p).
\end{split}
\end{equation}
\end{pro}
\pf  For completeness, we sketch the proofs.
First, it readily follows from (\ref{1.19})$_3$ and the  maximum principle that
\begin{equation}
\|\theta^{(0)}(t)\|_{L^\infty}\leq C,\quad\forall\ t\in[0,T].\label{3.6}
\end{equation}

Next, operating $A$ to both sides of (\ref{1.19})$_1$, multiplying it by $\uu^{(0)}$ in $L^2$ and integrating by parts, we have from (\ref{3.1}), (\ref{3.2}) and the Poincar${\rm\acute e}$'s and Cauchy-Schwarz's inequalities that
\begin{equation}\label{3.7}
\begin{split}
\|\nabla\uu^{(0)}\|_{L^2}^2\leq& C\left\|e^{-\tau A}\uu_0\right\|_{H^1}^2+C\left\|e^{-\tau A}\left(\kk \theta_{0}+\BB_0\cdot \nabla \BB_0+\partial_y\BB_0\right)\right\|_{L^2}^2\\
&+C\left(\| \theta^{(0)}\|_{L^2}^2+\|\BB^{(0)}\|_{L^2}^2\right)+\langle\BB^{(0)}\cdot\nabla \BB^{(0)},\uu^{(0)}\rangle\\
\leq&C\left(1+\| \theta^{(0)}\|_{L^2}^2+\|\BB^{(0)}\|_{L^2}^2\right)+\langle\BB^{(0)}\cdot\nabla \BB^{(0)},\uu^{(0)}\rangle.
\end{split}
\end{equation}

Similarly, multiplying (\ref{1.19})$_2$ and (\ref{1.19})$_3$ by $\BB^{(0)}$ and $\theta^{(0)}-(1-y)$ in $L^2$ respectively, we obtain after integrating by parts that
\begin{equation}\label{3.8}
\begin{split}
&\frac{d}{dt}\left(\|\BB^{(0)}\|_{L^2}^2+\|\theta^{(0)}\|_{L^2}^2\right)+\|\nabla\BB^{(0)}\|_{L^2}^2+\|\nabla\theta^{(0)}\|_{L^2}^2\\
&\quad\leq \frac{1}{2}\|\nabla\uu^{(0)}\|_{L^2}^2+C\left(1+\|\BB^{(0)}\|_{L^2}^2\right)+\langle\BB^{(0)}\cdot\nabla \uu^{(0)},\BB^{(0)}\rangle.
\end{split}
\end{equation}

Thus,  combining (\ref{3.7}) with (\ref{3.8}) and integrating by parts, we infer from the Gronwall's inequality that
\begin{equation}\label{3.9}
\sup_{0\leq t\leq T}\left(\|\BB^{(0)}\|_{L^2}^2+\|\theta^{(0)}\|_{L^2}^2\right)(t)+\int_0^T \|(\nabla\uu^{(0)},\nabla\BB^{(0)},\nabla\theta^{(0)})\|_{L^2}^2 \leq C.
\end{equation}

In a manner similar to the derivation of (\ref{2.7}), by (\ref{3.9}) we have
\begin{equation}\label{3.10}
\sup_{0\leq t\leq T}\|\BB^{(0)}(t)\|_{L^4}^4+\int_0^T\left(\||\BB^{(0)}||\nabla\BB^{(0)}|\|_{L^2}^2+\|\nabla|\BB^{(0)}|^2\|_{L^2}^2\right)dt\leq C,
\end{equation}
and hence, it is easily obtained  from (\ref{3.7}) that
\begin{equation}\label{3.11}
\|\uu^{(0)}(t)\|_{H^1}^2\leq C+C\|\BB^{(0)}(t)\|_{L^4}^4\leq C,\quad\forall\ t\in[0,T].
\end{equation}

Thanks to (\ref{3.2}) and (\ref{3.9}), we have
\begin{equation}\label{3.12}
\begin{split}
\|A\uu^{(0)}\|_{L^2} \leq& \|Ae^{-\tau A}\uu_0\|_{L^2}+\left\|e^{-\tau A}\mathbb{P}\left(\kk \theta_{0}+\BB_0\cdot \nabla \BB_0+\partial_ y\BB_0\right)\right\|_{L^2}\\
&+\left\|\mathbb{P}\left(\kk \theta^{(0)}+\BB^{(0)}\cdot\nabla \BB^{(0)}+\partial_y\BB^{(0)}\right)\right\|_{L^2}\\
\leq&C\left(1+ \||\BB^{(0)}||\nabla \BB^{(0)}|\|_{L^2}+\|\nabla\BB^{(0)}\|_{L^2}\right),
\end{split}
\end{equation}
so that, by (\ref{3.9}) and (\ref{3.10}) we deduce
\begin{equation}\label{3.13}
 \int_0^T\|\nabla^2\uu^{(0)}\|_{L^2}^2dt
\leq C+C\int_0^T\left(\||\BB^{(0)}||\nabla \BB^{(0)}|\|_{L^2}^2+\|\nabla\BB^{(0)}\|_{L^2}^2\right)dt\leq C,
\end{equation}
since it holds that $\|\nabla^2\uu^{(0)}\|_{L^2}=\|A\uu^{(0)}\|_{L^2}$ (see, for example, \cite{L1963,T2001}).

Similarly to  the proofs of (\ref{2.11}), (\ref{2.12}) and  (\ref{2.17}), by (\ref{3.9})--(\ref{3.13}) we can show that
\begin{equation}\label{3.14}
\sup_{0\leq t\leq T}\|(\nabla\BB^{(0)},\nabla\theta^{(0)})(t)\|_{{L^2}{\cap}{L^p}}+\int_0^T\left(\|(\BB^{(0)},\theta^{(0)})\|_{H^2}^2+\|(\BB^{(0)}_t,\theta^{(0)}_t)\|_{L^2}^2\right)dt\leq C,
\end{equation}
which, combined with (\ref{3.12}), yields
\begin{equation}\label{3.15}
\|\uu^{(0)}(t)\|_{H^2}\leq C+C\left(1+\|\BB^{(0)}(t)\|_{L^\infty}\right)\|\BB^{(0)}(t)\|_{H^1}\leq C,\quad\forall\ t\in[0,T].
\end{equation}

Therefore, collecting (\ref{3.9})--(\ref{3.15}) together finishes the proof of Proposition \ref{pro3.1}.\hfill$\square$

\subsection{Convergence from  $(\mathbf{u}^{(0)}, \mathbf{B}^{(0)},\theta^{(0)})$ to $(\mathbf{u}^{0}, \mathbf{B}^{0},\theta^{0})$}
In this subsection, we verify that as $\varepsilon\to0$, the solution of the effective dynamics (\ref{1.19})--(\ref{1.20}) is close to the one of the infinite Prandtl number model dynamics  (\ref{1.5})--(\ref{1.6}), based on the global (uniform) estimates stated in Propositions \ref{pro2.2} and \ref{pro3.1}. Indeed, it is  the infinite Prandtl number dynamics if the initial-layer corrections in (\ref{1.19})$_1$  are neglected.

To begin, noticing that
\begin{equation*}
\uu^0=A^{-1}\mathbb{P}\left(\kk \theta^{0}+\BB^{0}\cdot\nabla \BB^{0}+\frac{\partial\BB^{0}}{\partial y}\right),
\end{equation*}
we infer from (\ref{1.5}) and (\ref{1.19}) that
\begin{equation}\label{3.16}
\begin{split}
\uu^*\triangleq \uu^{(0)}-\uu^0=&e^{-\tau A}\uu_0-e^{-\tau A}A^{-1}\mathbb{P}\left(\kk \theta_{0}+\BB_0\cdot \nabla \BB_0+\BB_{0y}\right)\\
&+A^{-1}\mathbb{P}\left(\kk \theta^*+\BB^*\cdot\nabla \BB^{(0)}+\BB^0\cdot\nabla\BB^*+\BB^*_y\right),
\end{split}
\end{equation}
where $\BB^*\triangleq\BB^{(0)}-\BB^0$ and $\theta^*\triangleq\theta^{(0)}-\theta^0$ satisfy
\begin{equation}\label{3.17}
\frac{\partial \BB^*}{\partial t}+\uu^{(0)}\cdot \nabla\BB^*-\Delta\BB^*=-\uu^*\cdot \nabla \BB^0+\BB^*\cdot\nabla \uu^{(0)}+\BB^0\cdot \nabla\uu^*+\uu^*_y,
\end{equation}
and
\begin{equation}\label{3.18}
\frac{\partial\theta^*}{\partial t}+\uu^{(0)}\cdot \nabla \theta^*-\Delta \theta^*=-\uu^*\cdot \nabla \theta^0.
\end{equation}

Since it holds that $\BB^*\cdot\nabla \BB^{(0)}=\divg (\BB^*\otimes\BB^{(0)})$ and $\BB^0\cdot\nabla\BB^*=\divg(\BB^0\otimes\BB^*)$, using (\ref{3.2}), Propositions \ref{pro2.2} and \ref{pro3.1}, we easily deduce from (\ref{3.16}) that
\begin{equation}
\begin{split}
\|\uu^*\|_{H^1}
\leq&Ce^{-\tau}\left(\|\uu_0\|_{H^2}+\|\BB_0\|_{H^2}^2+\|\BB_0\|_{H^2}+\|\theta_0\|_{H^2}\right)\\
&+C\left(\|\theta^*\|_{L^2}+\|(\BB^{(0)},\BB^{0})\|_{L^\infty}\|\BB^*\|_{L^2}+\|\BB^*\|_{L^2}\right)\\
\leq&Ce^{-\tau}+C\left(\|\theta^*\|_{L^2}+\|\BB^*\|_{L^2}\right).
\end{split}\label{3.19}
\end{equation}

Multiplying (\ref{3.17}), (\ref{3.18}) by $\BB^*,\theta^*$ in $L^2$ respectively, integrating by parts, using Propositions \ref{pro2.2}, \ref{pro3.1} and the Poincar${\rm\acute{e}}$'s inequality, we have
\begin{equation*}
\begin{split}
&\frac{d}{dt}\|(\BB^*,\theta^*)\|_{L^2}^2+\|(\nabla\BB^*,\nabla\theta^*)\|_{L^2}^2\\
&\quad\leq C\|\uu^*\|_{L^4}\|(\nabla\BB^0,\nabla\theta^0)\|_{L^4}\|(\BB^*,\theta^*)\|_{L^2}+C\|\nabla\uu^*\|_{L^2}\|\BB^*\|_{L^2}\\
&\qquad+C\|\BB^*\|_{L^4}\|\nabla\uu^{(0)}\|_{L^4}\|\BB^*\|_{L^2}+C\|\BB^0\|_{L^\infty}\|\nabla\uu^*\|_{L^2}\|\BB^*\|_{L^2}\\
&\quad\leq C\|\uu^*\|_{H^1}\|(\BB^*,\theta^*)\|_{L^2}+C\|\nabla\BB^*\|_{L^2}\|\BB^*\|_{L^2},
\end{split}
\end{equation*}
which, combined with (\ref{3.19}) and the Cauchy-Schwarz's inequality, yields
\begin{equation}\label{3.20}
\frac{d}{dt}\|(\BB^*,\theta^*)\|_{L^2}^2+\|(\nabla\BB^*,\nabla\theta^*)\|_{L^2}^2\leq C_1\|(\BB^*,\theta^*)\|_{L^2}^2+C_2e^{-\tau}\|(\BB^*,\theta^*)\|_{L^2}.
\end{equation}

Keeping in mind that $\tau=t/\varepsilon$ and that $(\BB^*,\theta^*)|_{t=0}=0$, we deduce from (\ref{3.20}) that
\begin{equation}
\|(\BB^*,\theta^*)(t)\|_{L^2}\leq \frac{C_2\varepsilon}{2+C_1\varepsilon}e^{C_1t/2},\quad\forall\ t\in[0,T],\label{3.21}
\end{equation}
which, inserted  into (\ref{3.19}) and combined with (\ref{3.16}), shows that for any $0\leq t\leq T$,
\begin{equation}\label{3.22}
\|\mathbf{u}^{(0)}-\mathbf{u}^{0}-e^{-\tau A}\mathbf{u}_0+e^{-\tau A}A^{-1}\mathbb{P}\left(\mathbf{k} \theta_{0}+\mathbf{B}_0\cdot \nabla \mathbf{B}_0+\mathbf{B}_{0y}\right)\|_{H^1}
\leq C\varepsilon.
\end{equation}

To prove the $H^2$-convergence,  we first utilize (\ref{3.19}), (\ref{3.21}), Propositions \ref{pro2.2} and \ref{pro3.1} to deduce from (\ref{3.17}) and (\ref{3.18}) in a manner similar to the proof of (\ref{2.41}) that
\begin{equation*}
\begin{split}
\frac{d}{dt}\|(\nabla\BB^*,\nabla\theta^*)\|_{L^2}^2\leq& C\|(\nabla\BB^*,\nabla\theta^*)\|_{L^2}^2+C\|\uu^*\|_{H^1}^2\\
\leq& C\|(\nabla\BB^*,\nabla\theta^*)\|_{L^2}^2+C\varepsilon^2+Ce^{-2\tau},
\end{split}
\end{equation*}
from which we find
\begin{equation}\label{3.23}
\|(\nabla\BB^*,\nabla\theta^*)\|_{L^2}^2\leq C\varepsilon.
\end{equation}
and thus, it follows directly from (\ref{3.16}) that
\begin{equation}\label{3.24}
\begin{split}
&\|\mathbf{u}^{(0)}-\mathbf{u}^{0}-e^{-\tau A}\mathbf{u}_0+e^{-\tau A}A^{-1}\mathbb{P}\left(\mathbf{k} \theta_{0}+\mathbf{B}_0\cdot \nabla \mathbf{B}_0+\mathbf{B}_{0y}\right)\|_{H^2}^2\\
&\quad\leq\|A^{-1}\mathbb{P} (\kk \theta^*+\BB^*\cdot\nabla \BB^{(0)}+\BB^0\cdot\nabla\BB^*+\BB^*_y)\|_{H^2}^2\\
&\quad\leq C\left(\|\theta^*\|_{L^2}+\|\BB^*\|_{H^1}\right) \leq C\varepsilon.
\end{split}
\end{equation}

In short, collecting (\ref{3.21})--(\ref{3.24}) together, we arrive at
\begin{thm}
\label{thm3.1} For any given $T>0$, assume that $(\mathbf{u}^{0}, \mathbf{B}^{0},\theta^{0})$  and $(\mathbf{u}^{(0)}, \mathbf{B}^{(0)},\theta^{(0)})$  are the solutions of \eqref{1.5}--\eqref{1.6} and \eqref{1.19}--\eqref{1.20} on $\Omega \times[0, T]$, respectively. Then, there exists a positive constant $C$, independent of $\varepsilon$, such that for any $t\in[0,T]$,
\begin{equation}\label{3.25}
\|\mathbf{u}^{(0)}-\mathbf{u}^{0}-e^{-\tau A}\mathbf{u}_0+e^{-\tau A}A^{-1}\mathbb{P}\left(\mathbf{k} \theta_{0}+\mathbf{B}_0\cdot \nabla \mathbf{B}_0+\mathbf{B}_{0y}\right)\|_{H^1}
\leq C\varepsilon,
\end{equation}
\begin{equation}\label{3.26}
\|\mathbf{u}^{(0)}-\mathbf{u}^{0}-e^{-\tau A}\mathbf{u}_0+e^{-\tau A}A^{-1}\mathbb{P}\left(\mathbf{k} \theta_{0}+\mathbf{B}_0\cdot \nabla \mathbf{B}_0+\mathbf{B}_{0y}\right)\|_{H^2}
\leq C\varepsilon^{1/2},
\end{equation}
and
\begin{equation}\label{3.27}
\|(\mathbf{B}^{(0)}-\mathbf{B}^{0},\theta^{(0)}-\theta^0)\|_{L^2}\leq C\varepsilon,\quad \|\nabla(\mathbf{B}^{(0)}-\mathbf{B}^{0},\theta^{(0)}-\theta^0)\|_{L^2}
\leq C\varepsilon^{1/2}.
\end{equation}
\end{thm}

\subsection{Convergence from $(\mathbf{u}, \mathbf{B},\theta)$  to $(\mathbf{u}^{(0)}, \mathbf{B}^{(0)},\theta^{(0)})$}
This subsection aims to justify the limit from $(\uu,\BB,\theta)$ to $(\uu^{(0)},\BB^{(0)},\theta^{(0)})$, which is more complicated than the previous ones. To do this, let
$$
\tilde\uu=\uu-\uu^{(0)},\quad \tilde\BB=\BB-\BB^{(0)},\quad\tilde\theta=\theta-\theta^{(0)}.
$$

Then it is easily derived from (\ref{1.3})$_1$ and (\ref{1.19})$_1$ that
\begin{equation}\label{3.28}
\varepsilon\left(\frac{\partial\tilde\uu}{\partial t}+\uu\cdot\nabla\tilde\uu+\tilde\uu\cdot\nabla\uu^{(0)}\right)+\nabla p=\Delta\tilde\uu+\kk\tilde\theta+\BB\cdot\nabla\tilde\BB+\tilde\BB\cdot\nabla\BB^{(0)}+\frac{\partial \tilde\BB}{\partial y}+\ff,
\end{equation}
where $\ff$ is defined as follows:
\begin{equation}\label{3.29}
\begin{split}
\ff\triangleq&-\varepsilon\frac{\partial \uu^{(0)}}{\partial t}-\varepsilon \uu^{(0)}\cdot\nabla\uu^{(0)}+\Delta\uu^{(0)}+\kk\theta^{(0)}+\BB^{(0)}\cdot\nabla\BB^{(0)}+\frac{\partial\BB^{(0)}}{\partial y}\\
=&Ae^{-\tau A}\uu_0-e^{-\tau A}{\mathbb{P}}\left(\kk\theta_0+\BB_0\cdot\nabla\BB_0+\frac{\partial\BB_0}{\partial y}\right)\\
&-\varepsilon A^{-1}{\mathbb{P}}\frac{\partial}{\partial t}\left(\kk\theta^{(0)}+\BB^{(0)}\cdot\nabla\BB^{(0)}+\frac{\partial \BB^{(0)}}{\partial y}\right)\\
&-\varepsilon \uu^{(0)}\cdot\nabla\uu^{(0)}-A\uu^{(0)}+\kk\theta^{(0)}+\BB^{(0)}\cdot\nabla\BB^{(0)}+\frac{\partial\BB^{(0)}}{\partial y}+\nabla p^{(0)}\\
=&-\varepsilon A^{-1}{\mathbb{P}}\frac{\partial}{\partial t}\left(\kk\theta^{(0)}+\BB^{(0)}\cdot\nabla\BB^{(0)}+\frac{\partial \BB^{(0)}}{\partial y}\right)-\varepsilon \uu^{(0)}\cdot\nabla\uu^{(0)}+\nabla p^{(0)}.
\end{split}
\end{equation}

In the following, we prove the convergence of $\tilde\uu$.  First, multiplying (\ref{3.28}) by $\tilde\uu$ in $L^2$ and integrating by parts, we deduce
\begin{equation}\label{3.30}
\begin{split}
\frac{\varepsilon}{2}\frac{d}{dt}\|\tilde\uu\|_{L^2}^2+\|\nabla\tilde\uu\|_{L^2}^2=&-\varepsilon\langle\uu\cdot\nabla\uu^{(0)},\tilde\uu\rangle+
\langle\kk\tilde\theta+\BB\cdot\nabla\tilde\BB+\tilde\BB\cdot\nabla\BB^{(0)}+ \tilde\BB_{y},\tilde\uu\rangle\\
&-\varepsilon \langle A^{-1}{\mathbb{P}}\frac{\partial}{\partial t}(\kk\theta^{(0)}+\BB^{(0)}\cdot\nabla\BB^{(0)}+\BB^{(0)}_{y}),\tilde\uu \rangle\triangleq \sum_{i=1}^3I_i,
\end{split}
\end{equation}
where we have used the facts that $\divg\tilde\uu=0$, $\tilde\uu|_{y=0,1}=0$ and $\uu=\tilde\uu+\uu^{(0)}$.
Using Propositions \ref{pro2.1}, \ref{pro3.1} and the Poincar${\rm{\acute e}}$'s inequality, we have
\begin{equation*}
I_1\leq C\varepsilon\|\uu\|_{L^\infty}\|\nabla\uu^{(0)}\|_{L^2}\|\tilde\uu\|_{L^2}\leq C\varepsilon\|\nabla\tilde\uu\|_{L^2}
\leq\frac{1}{4}\|\nabla\tilde\uu\|_{L^2}^2+C\varepsilon^2,
\end{equation*}
and similarly,
\begin{equation*}
\begin{split}
I_2\leq&C\left(\|\tilde\theta\|_{L^2}\|\tilde\uu\|_{L^2}+\|(\BB,\BB^{(0)})\|_{L^\infty}\|\tilde\BB\|_{L^2}\|\nabla\tilde\uu\|_{L^2}
+\|\tilde\BB\|_{L^2}\|\nabla\tilde\uu\|_{L^2}\right)\\
\leq&\frac{1}{4}\|\nabla\tilde\uu\|_{L^2}^2+C\left(\|\tilde\BB\|_{L^2}^2+\|\tilde\theta\|_{L^2}^2\right).
\end{split}
\end{equation*}

In view of (\ref{3.4}) and Proposition \ref{pro3.1}, we obtain after integrating by parts that
\begin{equation*}
\begin{split}
I_3\leq& C\varepsilon\left(\|\theta^{(0)}_t\|_{H^{-1}}\|\tilde\uu\|_{L^2}+\|\BB^{(0)}\|_{L^\infty}\|\BB^{(0)}_t\|_{H^{-1}}\|\nabla\tilde\uu\|_{L^2}
+\|\BB^{(0)}_t\|_{H^{-1}}\|\nabla\tilde\uu\|_{L^2}\right)\\
\leq& \frac{1}{4}\|\nabla\tilde\uu\|_{L^2}^2+ C\varepsilon^2\left(\|\theta^{(0)}_t\|_{H^{-1}}^2+\|\BB^{(0)}_t\|_{H^{-1}}^2\right)\\
\leq& \frac{1}{4}\|\nabla\tilde\uu\|_{L^2}^2+ C\varepsilon^2,
\end{split}
\end{equation*}
where we have used (\ref{1.19}) and Proposition \ref{pro3.1}  to get that
\begin{equation*}
\|(\BB^{(0)}_t,\theta^{(0)}_t)\|_{H^{-1}}\leq \|(\BB^{(0)},\theta^{(0)})\|_{H^1}+\|\uu^{(0)}\|_{L^\infty}\|(\BB^{(0)},\theta^{(0)})\|_{L^2}+\|\uu^{(0)}\|_{L^2}\leq C.
\end{equation*}

Thus, substituting the estimates of $I_i$ ($i=1,2,3$) into (\ref{3.30}), we arrive at
\begin{equation}\label{3.31}
\begin{split}
\frac{\varepsilon}{2}\frac{d}{dt}\|\tilde\uu\|_{L^2}^2+\|\nabla\tilde\uu\|_{L^2}^2\leq C\left(\|\tilde\BB\|_{L^2}^2+\|\tilde\theta\|_{L^2}^2\right)+ C\varepsilon^2.
\end{split}
\end{equation}

Clearly, we need to deal with $\|(\tilde\BB,\tilde\theta)\|_{L^2}^2$. Indeed, by (\ref{1.3}) and (\ref{1.19}) we have
\begin{equation}\label{3.32}
\begin{dcases}
\frac{\partial \tilde\BB}{\partial t}+\uu\cdot \nabla \tilde\BB-\Delta\tilde\BB=-\tilde\uu\cdot \nabla \BB^{(0)}+\tilde\BB\cdot\nabla \uu+\BB^{(0)}\cdot \nabla\tilde\uu+\frac{\partial \tilde\uu}{\partial y},
\\
\frac{\partial \tilde\theta}{\partial t}+\uu\cdot \nabla\tilde\theta-\Delta\tilde\theta=-\tilde\uu\cdot \nabla \theta^{(0)},
\end{dcases}
\end{equation}
from which we easily get that
\begin{equation*}
\begin{split}
\frac{d}{dt}\|(\tilde\BB,\tilde\theta)\|_{L^2}^2+\|(\nabla\tilde\BB,\nabla\tilde\theta)\|_{L^2}^2\leq& C\left(1+\|(\uu,\BB^{(0)},\theta^{(0)})\|_{L^\infty}^2\right)\|(\tilde\uu,\tilde\BB)\|_{L^2}^2\\
\leq& C\|(\tilde\uu,\tilde\BB)\|_{L^2}^2,
\end{split}
\end{equation*}
where we have used  Propositions \ref{pro2.1} and \ref{pro3.1}. Thus,  for any $0\leq t\leq T$ one has
\begin{equation}\label{3.33}
\|(\tilde\BB,\tilde\theta)(t)\|_{L^2}^2+\int_0^t\|(\nabla\tilde\BB,\nabla\tilde\theta)\|_{L^2}^2ds\leq C\int_0^t\|\tilde\uu\|_{L^2}^2ds.
\end{equation}

In view of (\ref{3.31}), (\ref{3.33}) and the Poincar${\rm\acute{e}}$'s inequality, we have
\begin{equation*}
\frac{d}{dt}\left(e^{t/\varepsilon}\|\tilde\uu\|_{L^2}^2\right)\leq \varepsilon^{-1}e^{t/\varepsilon}\int_0^t\|\tilde\uu\|_{L^2}^2ds+C\varepsilon e^{t/\varepsilon},
\end{equation*}
which, integrated in time, yields (noting that $\tilde\uu|_{t=0}=0$)
\begin{equation*}
\|\tilde\uu(t)\|_{L^2}^2\leq C\int_0^t\|\tilde\uu\|_{L^2}^2ds+C\varepsilon^2,
\end{equation*}
and hence, by Gronwall's inequality we see that
\begin{equation}
\|\tilde\uu(t)\|_{L^2}^2\leq C\varepsilon^2,\quad \forall \ t\in[0,T].\label{3.34}
\end{equation}
This, combined with (\ref{3.31}) and (\ref{3.33}), also leads to
\begin{equation}\label{3.35}
\sup_{0\leq t\leq T}\|(\tilde\BB,\tilde\theta)(t)\|_{L^2}^2+\int_0^T\|(\nabla\tilde\uu,\nabla\tilde\BB,\nabla\tilde\theta)\|_{L^2}^2dt\leq C\varepsilon^2.
\end{equation}

Next, multiplying (\ref{3.32})$_1$ and (\ref{3.32})$_2$ by $\tilde\BB_t$ and $\tilde\theta_t$ in $L^2$ respectively, integrating by parts, using Propositions \ref{pro2.1}, \ref{pro3.1} and the Poincar${\rm\acute{e}}$'s inequality,  we obtain
\begin{equation*}
\begin{split}
&\frac{d}{dt}\left(\|\nabla\tilde\BB\|_{L^2}^2+\|\nabla\tilde\theta\|_{L^2}^2\right)+\left(\|\tilde\BB_t\|_{L^2}^2+\|\tilde\theta_t\|_{L^2}^2\right)\\
&\quad\leq C\|\uu\|_{L^\infty}^2\left(\|\nabla\tilde\BB\|_{L^2}^2+\|\nabla\tilde\theta\|_{L^2}^2\right)+C\|\nabla\tilde\uu\|_{L^2}^2+C\|\BB^{(0)}\|_{L^\infty}^2\|\nabla\tilde\uu\|_{L^2}^2\\
&\qquad+C\left(\|\nabla\BB^{(0)}\|_{L^4}^2+\|\nabla\theta^{(0)}\|_{L^4}^2\right)\|\tilde\uu\|_{L^4}^2+C\|\nabla\uu\|_{L^4}^2\|\tilde\BB\|_{L^4}^2\\
&\quad\leq C\left(\|\nabla\tilde\BB\|_{L^2}^2+\|\nabla\tilde\theta\|_{L^2}^2+\|\nabla\tilde\uu\|_{L^2}^2\right),
\end{split}
\end{equation*}
which, combined with (\ref{3.35}), results in
\begin{equation}\label{3.36}
\sup_{0\leq t\leq T}\left(\|\nabla\tilde\BB\|_{L^2}^2+\|\nabla\tilde\theta\|_{L^2}^2\right)+\int_0^T\left(\|\tilde\BB_t\|_{L^2}^2+\|\tilde\theta_t\|_{L^2}^2\right)dt\leq C\varepsilon^2,
\end{equation}
and moreover, it follows from  (\ref{3.32}) that
\begin{equation}\label{3.37}
\int_0^T\left(\|\tilde\BB\|_{H^2}^2+\|\tilde\theta\|_{H^2}^2\right)dt\leq C\varepsilon^2.
\end{equation}

Similarly, multiplying (\ref{3.28}) by $\tilde\uu_t$ in $L^2$ and integrating by parts, we obtain
\begin{equation*}
\begin{split}
\frac{d}{dt}\|\nabla\tilde\uu\|_{L^2}^2+\varepsilon\|\tilde\uu_t\|_{L^2}^2\leq& C\left(\varepsilon \|(\uu,\uu^{(0)})\|_{H^2}\|\tilde\uu\|_{H^1}+\varepsilon\|\uu^{(0)}\|_{H^2}^2+\|\tilde\theta\|_{L^2}+\|\tilde\BB\|_{H^1} \right)\|\tilde\uu_t\|_{L^2}\\
&+C\left(\|\BB\|_{L^\infty}\|\tilde\BB\|_{H^1}+\|\nabla\BB^{(0)}\|_{L^4}\|\tilde\BB\|_{H^1} \right)\|\tilde\uu_t\|_{L^2}\\
&+C\varepsilon\| A^{-1}{\mathbb{P}} (\kk\theta^{(0)}+\BB^{(0)}\cdot\nabla\BB^{(0)}+\BB^{(0)}_y)_t\|_{L^2}\|\tilde\uu_t\|_{L^2}\\
\leq& C\varepsilon\left(1+\|\theta^{(0)}_t\|_{H^{-1}}+\|\BB^{(0)}\|_{L^\infty}\|\BB^{(0)}_t\|_{L^2}+\|\BB^{(0)}_t\|_{L^2}\right)\|\tilde\uu_t\|_{L^2}\\
\leq& \frac{\varepsilon}{2}\|\tilde\uu_t\|_{L^2}^2+C\varepsilon\left(1+\|\BB^{(0)}_t\|_{L^2}^2\right),
\end{split}
\end{equation*}
where Propositions \ref{pro2.1}, \ref{pro3.1}, (\ref{3.4}) and  (\ref{3.34})--(\ref{3.36}) were used. Hence, by (\ref{3.37}) we find
\begin{equation}\label{3.38}
\sup_{0\leq t\leq T}\|\nabla\tilde\uu(t)\|_{L^2}^2+\varepsilon\int_0^T\|\tilde\uu_t\|_{L^2}^2dt\leq C\varepsilon+C\varepsilon\int_0^T\|\BB^{(0)}_t\|_{L^2}^2dt\leq C\varepsilon.
\end{equation}

Finally, it remains to derive the $H^2$-convergence of $\tilde\uu$, which is more subtle and needs more work. To do this, let
$$
\vv\triangleq \uu-\uu^{(0)}+A^{-1}{\mathbb{P}} \left(\kk\theta^{(0)}+\BB^{(0)}\cdot\nabla\BB^{(0)}+\BB^{(0)}_{y}\right).
$$

It is easy to deduce from (\ref{1.3})$_1$ and (\ref{1.19})$_1$ that
\begin{equation}
\varepsilon\vv_t+A\vv={\mathbb{P}} \left(\kk\theta +\BB \cdot\nabla\BB +\BB_{y}-\varepsilon \uu\cdot\nabla\uu\right).\label{3.39}
\end{equation}
Moreover, by direct calculations we have
\begin{equation*}
\begin{split}
\varepsilon \vv_t=&-A\uu+A\uu^{(0)}-{\mathbb{P}} \left(\kk\theta^{(0)}+\BB^{(0)}\cdot\nabla\BB^{(0)}+\BB^{(0)}_{y}\right)\\
&+{\mathbb{P}} \left(\kk\theta +\BB \cdot\nabla\BB +\BB_{y}-\varepsilon \uu\cdot\nabla\uu\right)\\
=&-A\uu+ Ae^{-\tau A}\uu_0-e^{-\tau A}{\mathbb{P}} \left(\kk\theta_0+\BB_0\cdot\nabla\BB_0+\BB_{0y}\right)\\
&+{\mathbb{P}}\left(\kk\theta +\BB \cdot\nabla\BB +\BB_{y}-\varepsilon  \uu\cdot\nabla\uu\right),
\end{split}
\end{equation*}
and hence,
\begin{equation}
\vv_t|_{t=0}=\lim_{t\to0}\vv_t= {\mathbb{P}}(\uu_0\cdot\nabla\uu_0)\in H^1.\label{3.40}
\end{equation}

Now, differentiating (\ref{3.39}) with respect to $t$ and multiplying it by $\vv_t$ in $L^2$, we deduce after integrating by parts that
\begin{equation}
\frac{\varepsilon}{2}\frac{d}{dt}\|\vv_t\|_{L^2}^2+\|\nabla\vv_t\|_{L^2}^2=\left\langle \left(\kk\theta +\BB \cdot\nabla\BB +\BB_{y}-\varepsilon \uu\cdot\nabla\uu\right)_t,\vv_t\right\rangle\triangleq I,\label{3.41}
\end{equation}
where the terms on the right-hand side can be estimated as follows, using Proposition \ref{pro2.1} and the Poincar${\rm\acute{e}}$'s inequality.
\begin{equation*}
\begin{split}
|I|
\leq& \left|\langle\kk\theta_t,\vv_t\rangle\right|+\left|\langle\BB_t\cdot\nabla\vv_t,\BB\rangle\right|+\left|\langle\BB\cdot\nabla\vv_t,\BB_t\rangle\right|\\
&+\left|\langle\BB_t,\vv_{ty}\rangle\right|+\varepsilon\left|\langle\uu_t\cdot\nabla\vv_t,\uu\rangle\right|
+\varepsilon\left|\langle\uu\cdot\nabla\vv_t,\uu_t\rangle\right|\\
\leq&\frac{1}{2}\|\nabla\vv_t\|_{L^2}^2+C\left(\|(\BB_t,\theta_t)\|_{L^2}^2+\varepsilon^2\|\uu_t\|_{L^2}^2\right)
\end{split}
\end{equation*}
Thus, by virtue of (\ref{3.40}) and Proposition \ref{pro2.1} we  get that
\begin{equation}
\varepsilon\sup_{0\leq t\leq T}\|\vv_t\|_{L^2}^2+\int_0^T\|\nabla\vv_t\|_{L^2}^2dt\leq C.\label{3.42}
\end{equation}

Using (\ref{2.24}), Proposition \ref{pro3.1} and (\ref{3.4}), we infer from (\ref{3.42})  that
\begin{equation}\label{3.43}
\begin{split}
\varepsilon\|\tilde\uu_t\|_{L^2}^2\leq& \varepsilon\|\vv_t\|_{L^2}^2+\varepsilon\|A^{-1}{\mathbb{P}} (\kk\theta^{(0)}+\BB^{(0)}\cdot\nabla\BB^{(0)}+\BB^{(0)}_{y})_t\|_{L^2}^2\\
\leq &C+C\varepsilon\left(\|\theta^{(0)}_t\|_{H^{-1}}+\|\BB^{(0)}\|_{L^\infty}\|\BB^{(0)}_t\|_{L^2}+\|\BB^{(0)}_t\|_{L^2}\right)\\
\leq &C+C\varepsilon\left(1+\|\BB^{(0)}_t\|_{L^2}\right).
\end{split}
\end{equation}

Similarly to the derivation of (\ref{2.21}), by Proposition \ref{pro3.1} we have
\begin{equation}\label{3.44}
\begin{split}
\frac{1}{2}\frac{d}{dt}\|\BB_t^{(0)}\|_{L^2}^2+\|\nabla\BB^{(0)}_t\|_{L^2}^2=&\langle\BB^{(0)}\cdot\nabla\uu^{(0)}_t,\BB^{(0)}_t\rangle +\langle\BB^{(0)}_t\cdot\nabla\uu^{(0)},\BB^{(0)}_t\rangle\\
&+\langle\uu^{(0)}_{yt},\BB^{(0)}_t\rangle
-\langle\uu^{(0)}_t\cdot\nabla\BB^{(0)},\BB^{(0)}_t\rangle\\
=&-\langle\BB^{(0)}\cdot\nabla\BB^{(0)}_t,\uu^{(0)}_t\rangle -\langle\BB^{(0)}_t\cdot\nabla\BB^{(0)}_t,\uu^{(0)}\rangle\\
&-\langle\uu^{(0)}_{t},\BB^{(0)}_{yt}\rangle
+\langle\uu^{(0)}_t\cdot\nabla\BB^{(0)}_t,\BB^{(0)}\rangle\\
\leq&\frac{1}{2}\|\nabla\BB^{(0)}_t\|_{L^2}^2+C\left(\|\uu^{(0)}_t\|_{L^2}^2+\|\BB^{(0)}_t\|_{L^2}^2\right)\\
\leq&\frac{1}{2}\|\nabla\BB^{(0)}_t\|_{L^2}^2+C\left(\|\tilde\uu_t\|_{L^2}^2+\|\uu_t\|_{L^2}^2+\|\BB^{(0)}_t\|_{L^2}^2\right).
\end{split}
\end{equation}

Inserting (\ref{3.43}) into (\ref{3.44}) and using (\ref{2.15}), by Gronwall's inequality we obtain
\begin{equation}\label{3.45}
\varepsilon\sup_{0\leq t\leq T}\|\BB_t^{(0)}(t)\|_{L^2}^2+\varepsilon\int_0^T \|\nabla\BB^{(0)}_t\|_{L^2}^2dt\leq
C+C\varepsilon\int_0^T\|\uu_t\|_{L^2}^2 dt\leq  C,
\end{equation}
which, together with (\ref{3.43}), immediately gives
\begin{equation}\label{3.46}
\varepsilon\sup_{0\leq t\leq T}\|\tilde\uu_t(t)\|_{L^2}^2\leq C
\end{equation}

Now, using Propositions \ref{pro2.1}, \ref{pro3.1}, (\ref{3.34})--(\ref{3.38}), (\ref{3.45}), (\ref{3.46}) and the estimates of the Stokes equations, we deduce from (\ref{3.28}) and (\ref{3.29}) that for any $0\leq t\leq T$,
\begin{equation}\label{3.47}
\begin{split}
\|\nabla^2\tilde\uu\|_{L^2}^2\leq &C\left(\|\tilde\theta\|_{L^2}^2+\|\BB\|_{L^\infty}^2\|\nabla\tilde\BB\|_{L^2}^2+\|\nabla\BB^{(0)}\|_{L^4}^2\|\tilde\BB\|_{L^4}^2
+\|\nabla\tilde\BB\|_{L^2}^2\right)\\ &+C\varepsilon^2\left(\|\tilde\uu_t\|_{L^2}^2+\|(\uu,\uu^{(0)})\|_{H^2}^2\|\tilde\uu\|_{H^1}^2+\|\uu^{(0)}\|_{H^2}^4\right)\\
&+C\varepsilon^2\|A^{-1}{\mathbb{P}} (\kk\theta^{(0)}+\BB^{(0)}\cdot\nabla\BB^{(0)}+\BB^{(0)}_{y})_t\|_{L^2}^2\\
\leq& C\varepsilon^2+C\varepsilon^2\left(1+\|\tilde\uu_t\|_{L^2}^2+\|\tilde\uu\|_{H^1}^2
+\|\theta_t^{(0)}\|_{H^{-1}}^2
+\|\BB_t^{(0)}\|_{L^2}^2\right)\\
\leq& C\varepsilon.
\end{split}
\end{equation}

Thus, collecting (\ref{3.34})--(\ref{3.36}), (\ref{3.38}) and (\ref{3.47}) together, we conclude that
\begin{thm}
\label{thm3.2} Let $(\mathbf{u}, \mathbf{B},\theta)$  and $(\mathbf{u}^{(0)}, \mathbf{B}^{(0)},\theta^{(0)})$  be the solutions of \eqref{1.3}--\eqref{1.4} and \eqref{1.19}--\eqref{1.20} on $\Omega \times[0, T]$ with $0<T<\infty$, respectively. Then, there exists a positive constant $C$, independent of $\varepsilon$, such that for any $t\in[0,T]$,
\begin{equation}\label{3.48}
\|\mathbf{u}-\mathbf{u}^{(0)}\|_{L^2}+\|(\mathbf{B}-\mathbf{B}^{(0)},\theta-\theta^{(0)})\|_{H^1}\leq C\varepsilon
\end{equation}
and
\begin{equation}\label{3.49}
\|\nabla(\mathbf{u}-\mathbf{u}^{(0)})\|_{L^2}+\|\nabla^2(\mathbf{u}-\mathbf{u}^{(0)})\|_{L^2}\leq C\varepsilon^{1/2}.
\end{equation}
\end{thm}

\noindent{\it Proof of Theorem \ref{thm1.2}.}  Now, combining Theorems \ref{thm3.1} and \ref{3.2}, we readily obtain the desired convergence results stated in (\ref{1.21})--(\ref{1.23}) of Theorem \ref{thm1.2}.\hfill$\square$

\vskip 4mm

\noindent{\bf Acknowledgment.} This work is partially
supported by NNSFC (Grant No. 11271306), the Natural Science
Foundation of Fujian Province of China (Grant No. 2015J01023), and the Fundamental Research Funds for the
Central Universities (Grant No. 20720160012).


\vskip 4mm


 {\small
\begin{thebibliography}{99}

\bibitem{Br2014}
Y. Brenier, Topology-preserving diffusion of divergence-free vector fields and magnetic relaxation. Comm. Math. Phys. 330 (2014), no. 2, 757-770.

\bibitem{Bu1989}
F. H. Busse, Fundamentals of Thermal Convection. Mantel convection: plate tectonics and fluid dynamics, 23-95. W. R. Peltier, ed. The Fluid Mechanics of Astrophysics and Geophysics, vol 4. Gordon and Breach, New York, 1989.

\bibitem{CW2011}
C. S. Cao, J. H. Wu, Global regularity for the 2D MHD equations with mixed partial dissipation and magnetic diffusion. Adv. Math. 226 (2011), no. 2, 1803-1822.

\bibitem{CWY2014}
C. S. Cao, J. H. Wu, B. Q. Yuan, The 2D incompressible magnetohydrodynamics equations with only magnetic diffusion. SIAM J. Math. Anal. 46 (2014), no. 1, 588-602.


\bibitem{C1952}
S. Chandrasekhar, On the inhibition of convection by a magnetic field. Phil. Mag. Ser. 7, 43 (1952),  501-532.

\bibitem{C1954}
S. Chandrasekhar, On the inhibition of convection by a magnetic field. II. Phil. Mag. Ser. 7, 45 (1954), 1177-1191.

\bibitem{C1961}
S. Chandrasekhar, Hydrodynamic and Hydromagnetic Stability. Clarendon, Oxford, 1961.

\bibitem{CD1999}
P. Constantin, C. R. Doering, Infinite Prandtl number convection. J. Statist. Phys. 94 (1999), no. 1-2, 159-172.

\bibitem{CHP2001}
P. Constantin, C. Hallstrom, V. Poutkaradze, Logarithmic bounds for infinite Prandtl number rotating convection. J. Math. Phys. 42 (2001), 773-783.

\bibitem{Da2001}
P. A. Davidson, An Introduction to Magnetohydrodynamics, Combridge University Press, 2001.

\bibitem{DC2001}
C. R. Doering, P. Constantin, On upper bounds for infinite Prandtl number convection with or without rotation. J. Math. Phys. 42 (2001),  784-795.

\bibitem{Ga1985}
G. P. Galdi, Nonlinear stability of the magnetic BB\'enard problem via a generalized energy method. Arch. Rational Mech. Anal. 87 (1985), no. 2, 167-186.

\bibitem{G1998}
A. V. Getling, Rayleigh-B\'enard Convection. Structure and Dynamics. Advanced Series in Nonlinear Dynamics, 11. World Scientific, River Edge, N.J. 1998.

\bibitem{GL2000}
S. Grossmann, D. Lohse, Scaling in thermal convection: a unifying theory. J. Fluid Mech. 407 (2000), 27-56.

\bibitem{Ho1995}
M.H. Holmes, Introduction to Perturbation Methods. Springer, New York, 1995.

\bibitem{L1963}
O. A. Ladyzhenskaya, The mathematical theory of viscous incompressible flow.  Silverman Gordon and Breach Science Publishers, New York-London 1963.

\bibitem{Ma2003}
A. Majda, Introduction to PDEs and Waves for the Atmosphere and Ocean. Courant Lecture Notes in Mathematics. AMS Providence, R. I., 2003.

\bibitem{MRR2014}
D. J. Mccormick, J. C. Robinson, J. L. Rodrigo, Existence and uniqueness for a coupled parabolic-elliptic model with applications to magnetic relaxation. Arch. Rational Mech. Anal. 214 (2014), 503-523.

\bibitem{MR1994}
G. Mulone, S. Rionero, On the stability of the rotating B\'enard problem. Bull. Tech. Univ. Istanbul 47 (1994), 181-202.

\bibitem{MR2003}
G. Mulone, S. Rionero, Necessary and sufficient conditions for nonlinear stability in the magnetic BšŠnard problem. Arch. Ration. Mech. Anal. 166 (2003), no. 3, 197-218.

\bibitem{N1955}
Y. Nakagawa, An experiment on the inhibition of thermal convection by a magnetic field. Nature, 175 (1955), 417-419.

\bibitem{N1957}
Y. Nakagawa, Experiments on the inhibition of thermal convection by a magnetic field. Proc. Royal Soc. London A, 240 (1957), 108-113.

\bibitem{Pa1979}
E. N. Parker, Cosmical Magnetic Fields, Clarendon, Oxford, 1979.

\bibitem{R1997}
S. Rionero, On magnetohydrodynamic stability, Quaderni di Matematica, 1 (1997), 347-376.

\bibitem{ST1983}
M. Sermange, R. Temam, Some mathematical questions related to the MHD equations. Commun. Pure Appl. Math. 36 (1983), 635-64.

\bibitem{S1994}
E. D. Siggia, High Rayleigh number convection. Annual Review of Fluid Mechanics, Vol. 26, 137-168. Annual Reviews, Palo Alto, Calif., 1994.

\bibitem{T2001}
R. Temam, Navier-Stokes equations. Theory and Numerical Analysis. AMS Providence, R. I., 2001.

\bibitem{T1951}
W. B. Thompson, Thermal convection in a magnetic field. Phil. Mag. Ser. 7, 42 (1951), 1417-1432 .



\bibitem{Wang2004}
X. M. Wang, Infinite Prandtl number limit of Rayleigh-Benard convection. Comm. Pure Appl. Math. 57 (2004), no. 10, 1265-1282.

\bibitem{Wang2007}
X. M. Wang, Asymptotic behavior of the global attractors to the Boussinesq system for Rayleigh-Benard convection at large Prandtl number. Comm. Pure Appl. Math. 60 (2007), no. 9, 1293-1318.

\end{thebibliography}}
\end{document}